\documentclass{article}

\usepackage{authblk}

\usepackage{natbib}
\usepackage{amssymb}
\usepackage{amsmath}
\usepackage{amsfonts}
\usepackage{amsthm}
\newtheorem{remark}{Remark}

\usepackage{graphicx}

\title{A short perspective on 
a posteriori error control and adaptive discretizations} 

\author[1]{Roland Becker}%
\author[2]{St\'ephane P.A. Bordas}
\author[3]{Franz Chouly}
\author[4,5]{Pascal Omnes}

\affil[1]{LMAP, University of Pau, IPRA BP 1155, Av. de l’Université, 64013 Pau, France}%

\affil[2]{Institute for Computational Engineering, Faculty of Science, Technology and Communication, University of Luxembourg, Luxembourg}%

\affil[3]{University of the Republic, Faculty of Science, Center of Mathematics, 11400 Montevideo, Uruguay}

\affil[4]{Université Sorbonne Paris Nord, LAGA, CNRS UMR 7539, Institut Galilée, 99 Av. J.-B. Clément, 93430, Villetaneuse, France}%
\affil[5]{Université Paris-Saclay, CEA, Service de Génie Logiciel pour la Simulation, 91191, Gif-sur-Yvette, France}

\date{\today}

\begin{document}

\maketitle


\begin{abstract}
   Error control by means of a posteriori error estimators or indicators and adaptive discretizations, such as adaptive mesh refinement, have emerged in the late seventies. Since then, numerous theoretical developments and improvements have been made, as well as the first attempts to introduce them into real-life industrial applications. The present introductory chapter provides an overview of the subject, highlights some of the achievements to date and discusses possible perspectives.
\medskip

\textit{This preprint corresponds to the Chapter 1 of volume 58 in AAMS, Advances in Applied Mechanics (to appear).}

\medskip

   
\end{abstract}



\section{Introduction}

The numerical solution of partial differential equations 
is a crucial aspect of computational science and engineering, finding applications in diverse fields of physics and engineering.
However, achieving accurate and efficient solutions to partial differential equations is often challenging and a posteriori error estimators are crucial tools that evaluate the accuracy and enhance computational efficiency of numerical methods. Unlike a priori error estimators, which are not computable, these estimators assess solutions after computation, providing a pragmatic evaluation of accuracy.

A posteriori estimators work  as evaluative tools, systematically analysing error distributions to offer insight into numerical approximation fidelity. This shift from anticipatory to retrospective assessment enables practitioners to iteratively refine solutions to align them more closely with the intricate details of the underlying mathematical models.

This article explores the conceptual landscape of a posteriori error estimators, examining their theoretical foundations, applications, and consequential impact on partial differential equations problem-solving. The discussion invites an exploration of numerical accuracy, navigating the terrain where mathematical abstraction converges with computational reality.

\section{A few basic notions}

Before entering the core of the topic, let us describe a few useful basic notions, in a very general setting and without entering into details. For reviews about a posteriori error estimation, see, e.g., \cite{eriksson1985,carstensen2010,NochettoVeeser12,chamoin2023}.

\subsection{Numerical approximation of a mathematical model}

For several decades, computers have been used in routine calculations to provide an approximate solution to the sophisticated mathematical models in current engineering practice, which are, for example, (linear or non-linear) partial differential equations  supplemented by (linear or non-linear) boundary conditions and by initial conditions, and possibly coupled to other equations, and now, to data (even in real-time).

Let us call $u$ the exact solution of this mathematical model, which belongs to a set $S$ of admissible solutions.
There are now many techniques for transforming the original equation through a sophisticated pipeline, so that, in the end, a concrete approximation to the solution $u$ is provided by a computer in the form of a collection of numbers $\left( U_1,\ldots,U_N \right) \in \mathbb{R}^N$ which can be used to recover a function $u_{\mathrm{computer}}(U_1,\ldots,U_N)$ (hopefully) close enough to $u$.

\begin{figure}[!ht]		
    \begin{center}	
    \includegraphics[width=\textwidth,keepaspectratio=true]{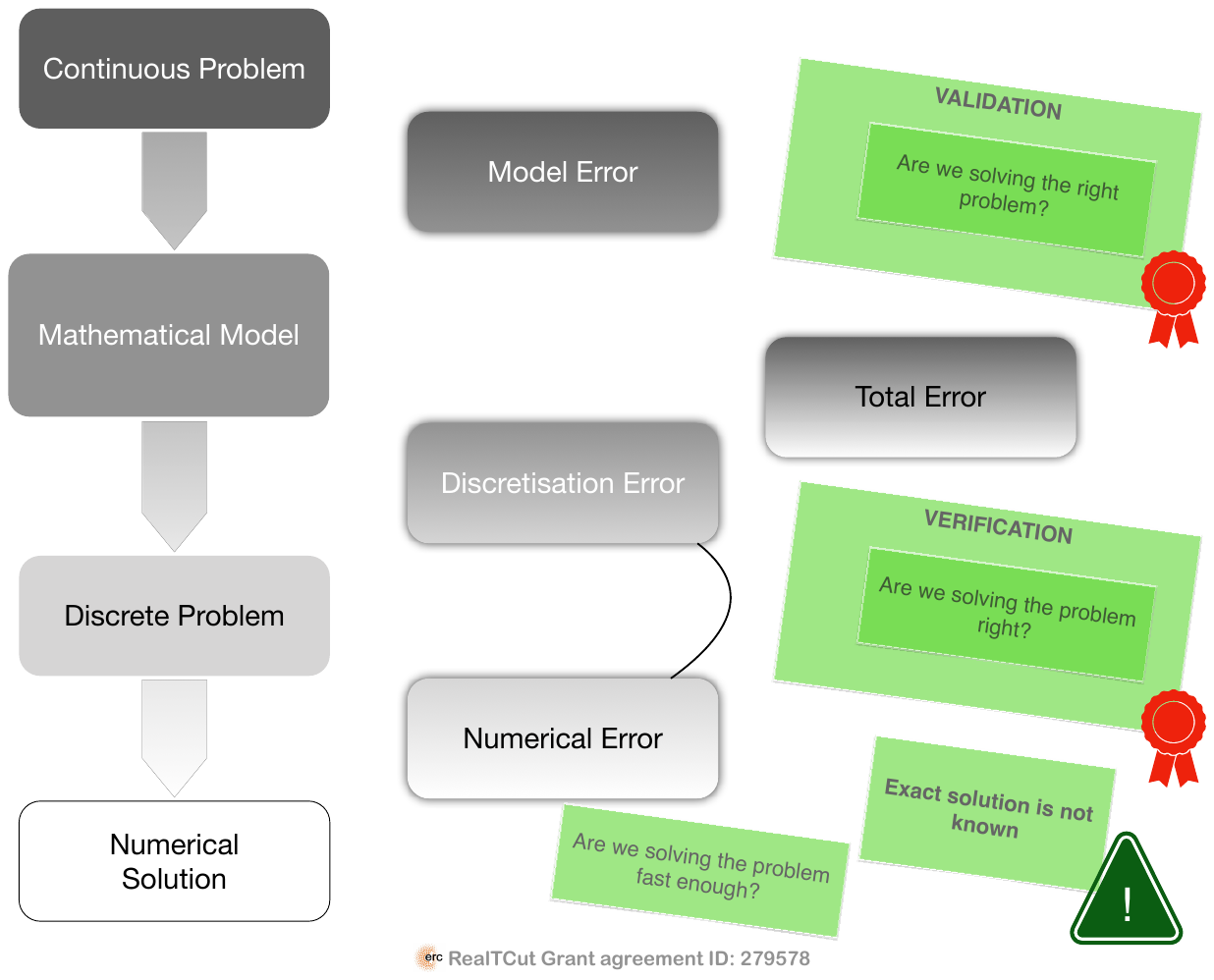}
    \end{center}	
    \caption{Mathematical modelling and sources of error. This figure showcases the steps involved in mathematical modelling, which is a process used to represent real-world phenomena using mathematical equations. The diagram highlights the different stages of mathematical modelling, including problem formulation, model development, parameter estimation, model validation, and prediction.
Additionally, the figure identifies potential sources of error that may impact the accuracy and reliability of the mathematical model. These sources of error can arise from various factors, such as measurement uncertainties, assumptions made during model development, limitations in data availability, simplifications in model assumptions, and uncertainties in parameter estimation.
The figure serves as a visual representation of the complexities and challenges associated with mathematical modelling, emphasising the need for careful consideration of potential sources of errors to ensure the robustness and validity of the model's results. It underscores the importance of thorough validation and verification processes to enhance the accuracy and reliability of mathematical models, which are crucial for decision-making, prediction, and understanding complex systems in various fields of science, engineering, and beyond.}
\label{fig:mathModelling}
\end{figure}

\begin{figure}[!ht]		
    \begin{center}	
    \includegraphics[width=\textwidth,keepaspectratio=true]{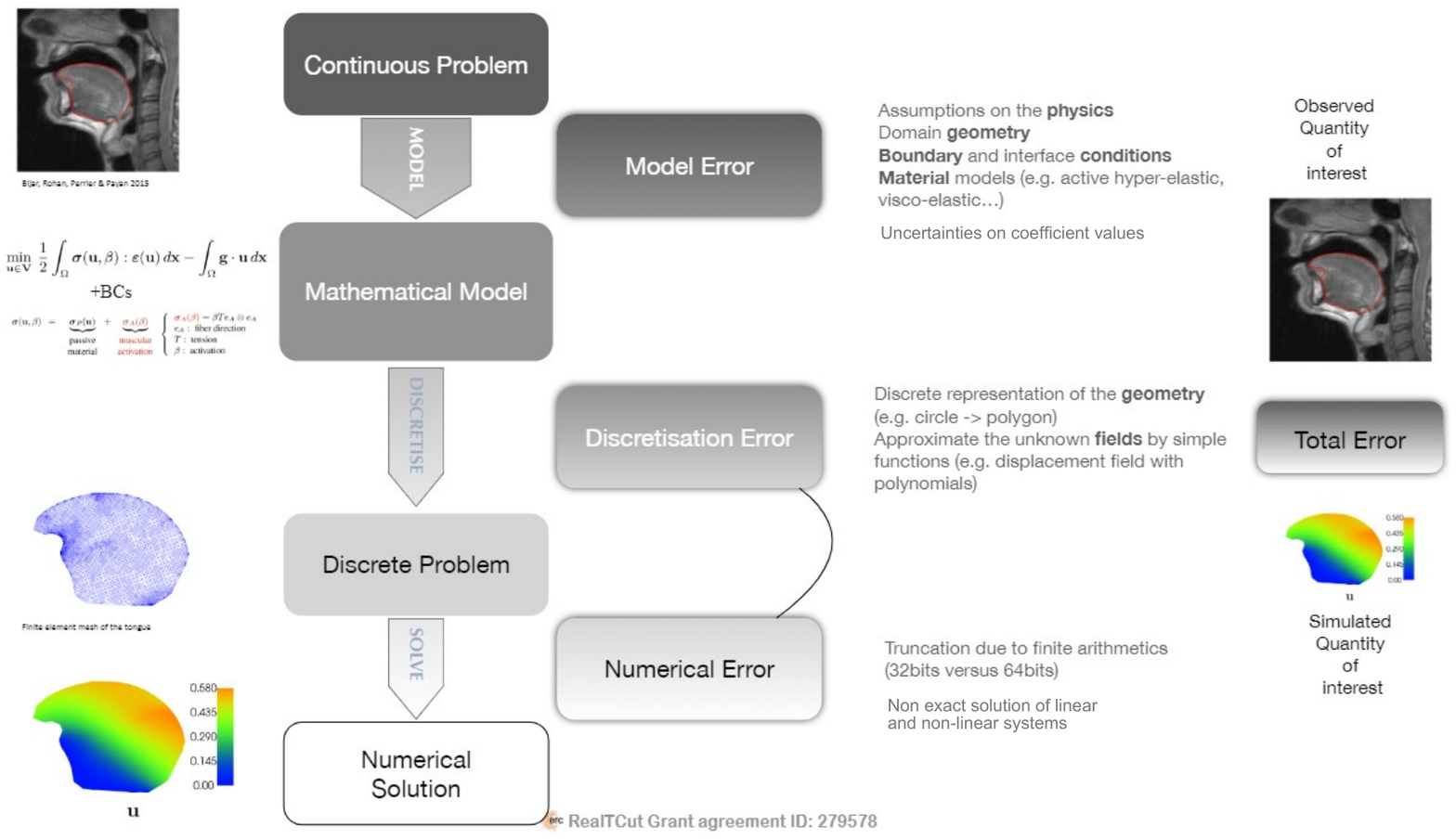}
    \end{center}	
    \caption{This figure showcases the process of mathematical modelling along with the identification and estimation of potential sources of error. The diagram illustrates the steps involved in mathematical modelling, including problem formulation, model development, parameter estimation, model validation, and prediction.
Furthermore, the figure highlights the importance of error estimation in the modelling process. It identifies potential sources of error, such as measurement uncertainties, assumptions made during model development, limitations in data availability, simplifications in model assumptions, and uncertainties in parameter estimation. The figure emphasises the need to account for and quantify these sources of error in order to assess the reliability and accuracy of the mathematical model.
The figure serves as a visual representation of the comprehensive approach to mathematical modelling, which includes not only model development but also thorough error estimation to enhance the robustness and validity of the model's results. It underscores the significance of error estimation in improving the quality of mathematical models and their applicability in various fields of science, engineering, and beyond.} 
\label{fig:mathModelling_2}
\end{figure}

Figure \ref{fig:mathModelling} and Figure \ref{fig:mathModelling_2} below provide an overview of the process of mathematical modelling. 

\subsection{Galerkin methods and the discretization error}
\label{sub:galerkin}

To be illustrative, we consider (Petrov-)Galerkin methods, that include a broad class of numerical approximation techniques based on the weak form of the mathematical model. This includes particularly the case of well-known Finite Element Methods (Lagrange, mixed, etc), see, e.g., 
\cite{brenner2008,ciarlet2002,ern2021a,quarteroni1994,szabo2021}
but also many other new methods, among which:
\begin{itemize}
    \item  discontinuous Galerkin (DG) methods, 
    see for instance the pioneering works \cite{arnold-82,lesaint1974} and the recent monograph 
    \cite{dipietro2012};
\item 
recent polytopal 
methods
such as 
Conforming Polygonal Finite Elements, see, e.g., \cite{sukumar2004};
Hybrid Discontinuous Galerkin (HDG), see, e.g., \cite{ern2016}, Hybrid High Order (HHO) methods, see, e.g.,  \cite{cicuttin2021,ern2016,dipietro2020,lemaire2021}, 
the Weak Galerkin Method, see, e.g., \cite{dong2022},
the Virtual Element Method (VEM), see, e.g. \cite{beirao2013vem,lemaire2021};
the Smooth Finite Element Method (SFEM), 
see, e.g., \cite{liu2007,nguyen-xuan-2008};
modified discontinuous Galerkin with static condensation, see, e.g., \cite{lozinski2019poly};


\item
Discontinuous Petrov Galerkin (DPG) methods, 
see, e.g., 
\cite{demkowicz2010,demkowicz2012};

\item 
IsoGeometric Analysis (IGA) and variants, see, e.g., \cite{atroshchenko2018,cottrell2009,nguyen2015isogeometric}, motivated by the link between computer aided design and numerical simulation;

\item unfitted finite elements or geometrically nonconforming finite elements, where the mesh boundary and the domain boundary do not match, as it occurs in fictitious domain methods, the eXtended Finite Element Method (XFEM) or the cut Finite Element Method (cutFEM): see, {e.g.}, \cite{bordas2023partition,burman2015,duprez2020fifem,glowinski-pan-periaux-94,haslinger2009,moes2006,nguyen2008meshless,peskin2002};

\item reduced basis techniques, such as Reduced Order modelling (ROM) or Proper Orthogonal Decomposition (POD), see, e.g., \cite{grepl2007,hesthaven2016,kerfriden2011};

\item spectral methods, see, e.g., \cite{bernardi1997} or \cite{canuto2006};

\item wavelet-based discretization, see, e.g., \cite{bertoluzza1994,bertoluzza1995,cohen2000,cohen2001,monasse1998}. 

\end{itemize}

Of course, the above list is far from exhaustive and the discussion can be extended, with appropriate modifications, to other classes of approximation techniques such as finite differences, finite volumes, collocation methods, etc, which are not based on a variational paradigm.

So, for a Petrov-Galerkin method, let us suppose the exact solution $u$ satisfies (for instance) the weak problem
\begin{equation}
    \label{weak}
    u \in V:\ 
    a(u;v)=L(v) \qquad \forall v\in W,
\end{equation}
where $V, W$ are some function spaces 
and $L$ is a linear form.
The form 
\[
a : V \times W \rightarrow \mathbb{R}
\]
is possibly nonlinear in the first variable.

\begin{remark}
The above formalism \eqref{weak} does not encompass nonlinear nonsmooth problems related to variational inequalities for instance, such as contact, friction, plasticity. However this is not critical for the present discussion, and to see how extensions for variational inequalities can be carried out, see, e.g., \cite{han2005book} or some contributions in these volumes, among which \cite{bartels2024exact,gustafsson2024adaptive,repin2024posteriori}.

. 
\end{remark}

A discrete Petrov-Galerkin method consists in approximating Problem~\eqref{weak} by a simpler problem in finite dimensional vector spaces:
\begin{equation}
    \label{weakRG}
    u_N \in V_N:\ 
    a_{G}(u_N;v_N)=L_{G}(v_N) \qquad \forall v_N\in W_N,
\end{equation}
where $V_N$ and $W_N$ are \emph{finite dimensional spaces} of dimension $N$, and $a_{G}$, resp. $L_{G}$, is a form that mimics $a$, resp. $L$ (in the simplest situations, we can take $a_G = a$, resp. $L_G=L$).

\begin{remark}
For many methods, the trial space $V_N$ and the test space $W_N$ are the same ($V_N=W_N$) and this corresponds to standard Galerkin (or standard Ritz-Galerkin) methods. The case where the spaces differ is often associated with the terminology of Petrov-Galerkin methods, or non-standard (Ritz-)Galerkin methods, see \cite{ern-guermond-04}.
For a historical perspective about this class of methods, see \cite{gander2012} and references therein, which emphasise W. Ritz's outstanding contributions related to the approximation of the spectral biharmonic problem for Chladni figures.    
\end{remark}
 As a result, $u_N$ can be represented in a basis of $N$ functions, and a computer can provide the $N$ values of the components when Problem \eqref{weakRG} is solved.

The challenge now is to design the spaces $V_N$ and $W_N$ carefully enough to calculate a \emph{discrete solution}~$u_N$ that represents the exact solution $u$ as accurately as possible. To measure the difference between $u$ and $u_N$, the concept of discretisation error $\epsilon_N$  is introduced.
A natural and simple definition can be
\begin{equation}
\label{def:discretizationerror}
\epsilon_N(u,u_N) := \| u - u_N \|_{V},
\end{equation}
where $\| \cdot \|_V$ is the natural norm associated with the Banach or Hilbert structure of~$V$ (Sobolev norm or energy norm). Of course, other possibilities exist motivated by practical applications, see paragraph~\ref{goal} below.

Since $N$ is linked to the available computing resources, a practitioner may wish to achieve the lowest possible value for the discretisation error $\epsilon_N$, while keeping~$N$ as small as possible. 
For example, and ideally, one would like $\epsilon_N$ to be 
of the order of 
machine precision ($\epsilon_N \simeq 10^{-16}$ for double precision arithmetic), and $N$ to be such that the solution $u_N$ is obtained in real time (a few milliseconds for example). Today, even with the enormous progress made in computing power, this objective remains a challenge, particularly for industrial applications where the geometry and mathematical model can be very complex. In general, a compromise has to be found between acceptable accuracy and acceptable use of computing resources, and this compromise is highly dependent on the context and the targeted applications.
To achieve acceptable accuracy, 
the main problem one encounters is that $\epsilon_N$ is unknown, because the exact solution $u$ itself is unknown. 

Another issue is that the solution $u_N$ to Problem~\eqref{weakRG} is not exactly the solution 
$u_{\textrm{computer}}(U_1,\ldots,U_N)$ delivered by a computer. Indeed, the code written to provide $u_{\textrm{computer}}(U_1,\ldots,U_N)$ also contains other approximations, that are mostly due to: numerical integration, iterative solvers for nonlinear and/or linear systems, and finite precision arithmetic, see, e.g., \cite{ciarlet2002,ern-guermond-04} in the context of finite element methods. This fact introduces another layer of errors, \emph{numerical errors}, that can be quantified as, for instance:
\[
\| u_{\textrm{computer}}(U_1,\ldots,U_N) - u_N \|_D
\]
where $\| \cdot \|_D$ is a convenient norm to measure this in finite dimensions (it does not have to be necessarily the same norm as for the discretization error).

In general, it is assumed that the numerical errors are of very small magnitude in comparison to discretization errors. However, in some situations they need to be taken into account: for instance a gradient conjugate solver can be stopped after a few iterations to save computation time. See for instance \cite{ern2013} for a technique that can resolve this issue.

Finally we can get back to Figure \ref{fig:mathModelling} and Figure \ref{fig:mathModelling_2} in which the discretization error and numerical error are depicted, and complemented with the model error that encompasses all the discrepancies between the actual physical system that needs to be modelled and the idealised mathematical model (this important topic is far beyond the scope of this short overview).
%

\subsection{Error estimators and error control}

Since  the late 1970s, pioneering work within the finite element community has shown how it is possible to compute a quantity $\eta_N(=\eta_N(u_N))$ that depends only of the discrete solution $u_N$ and 
that is equivalent to
the discretization error:
\begin{equation}
\label{eta}
    \eta_N(u_N) \simeq 
    \epsilon_N(u,u_N).
\end{equation}
Classical references are \cite{babuska1978}, \cite{ladeveze1983}, \cite{eriksson1985}, \cite{bank1985}, and the books \cite{ainsworth2000} and \cite{verfurth2013}.
For simple mathematical models, we can take
\begin{equation}
    \label{dualnorm}
    \eta_N(u_N)=
    \| L_G(\cdot) - a_G(u_N;\cdot) \|_\star
\end{equation}
where $ L_G(\cdot) - a_G(u_N;\cdot)$ is the residual associated with Problem~\eqref{weakRG} and $\|\cdot\|_\star$ is a dual norm.
A problem with \eqref{dualnorm}
comes from the norm $\|\cdot\|_\star$, that is not computable. 
A highly desirable property of an estimator is not only its computability from $u_N$, but it should allow for an implementation, such that the 
necessary time for computation is negligible in comparison to the time needed to solve Problem~\eqref{weakRG}.

An important mathematical property to have a trustworthy estimator is its \emph{reliability}: if it can not provide the exact discretization error, it should at least provides an \emph{upper bound} of the discretization error, in the sense that, there exists $C>0$ independent of $N$ and $u_N$, such that
\begin{equation}
    \label{upperbound}
    \epsilon_N(u,u_N) \leq C\, \eta_N,
\end{equation}
for all $u \in V$ and $u_N \in V_N$. 
See for instance \cite{verfurth1999,veeser2009,veeser2012} for a discussion and details in the case of the residual error estimate.
If the constant $C$ in \eqref{upperbound} is equal to $1$, we have an \emph{exact guaranteed upper bound}:
\begin{equation}
    \label{guaranteedupperbound}
    \epsilon_N(u,u_N) \leq \eta_N,
\end{equation}
which is of high practical interest since it ensures the discrete solution $u_N$ is approximated with a degree of accuracy that is known. See for instance { \cite{neitt2004,Vohralik2007SIAMJNA,braess2008} 
and a discussion in \cite{verfurth2009} about these ``constant-free'' estimates.}
This paves the way for \emph{certified numerical methods} in which one can ensure that the discretization error is below a known threshold.
Indeed, note that if \eqref{upperbound} is satisfied only, the estimator $\eta_N$ can underestimate the error for instance. 

Last but not least, it is usual to define the effectivity index of an estimator as
\begin{equation}
\label{efficiency}
    \mathrm{eff} := \frac{\eta_N(u_N)}{\epsilon_N(u,u_N)}.
\end{equation}
The above definition is motivated by \eqref{eta} and is meant to 
describe the accuracy of an estimator. Hopefully, it should not vary too much with $N$ and $u$, and should be close to 1. 
First, if (\ref{upperbound}) holds, this implies that
\begin{equation}
    \mathrm{eff} \geq \frac1{C}.
\end{equation}
Clearly, the boundedness of the effectivity index is equivalent to the lower bound
\begin{equation}
    \label{lowerbound}
    \eta_N \leq c\, \epsilon_N(u,u_N)
\end{equation}
for some $c>0$ independent of $u$ and $N$; we then have $\textrm{eff}\le c$.
%
Lower bounds are often called effectivity and have been an important research topic, see \cite{verfurth2013}. In the context of the finite element method, they generally even have a local form, see {\bf (\ref{efficiencyloc}) 
}.

Asymptotic exactness of an estimator refers to
\begin{equation}
    \label{assymptotiexact}
    \lim_{N\to\infty} \mathrm{eff}_N  = 1.
\end{equation}
This property depends on $u$, the method to compute $u_N$ and the sequence of spaces~$V_N$.
In special cases, estimators such as the Bank and Weiser type, are 
known to yield (\ref{assymptotiexact}), see \cite{babuska1992,bank1985,duran1991} for details.

\subsection{Error estimators and mesh refinement}
In addition to assess the accuracy of a simulation, error estimators are used to iteratively improve the numerical approximation. The typical form of the loop is
\begin{center}
    \textbf{Solve} $\quad\rightarrow\quad$ \textbf{Estimate} $\quad\rightarrow\quad$ 
    \textbf{Mark} $\quad\rightarrow\quad$  
    \textbf{Adapt} $\quad\rightarrow\quad\cdots$
\end{center}
where \textbf{Solve} refers to solution of the discrete problem for given approximation spaces, see, e.g., Problem \eqref{weakRG}, and the estimator is then used to adapt the approximation spaces. The algorithm produces a sequence of {meshes $\mathcal K_{\ell}$ 
and associated} spaces
$V_{\ell}$, solutions $u_{\ell}$ and estimators $\eta_{\ell}$.
The global upper bound justifies the use of the estimator $\eta_\ell$ as a stopping criterion.

In the finite element method, the approximation spaces are built from meshes and the \textbf{Adapt} step usually consists in local mesh-refinement.
Since the spatial domain associated with the mathematical problem is subdivided into cells $K$, say triangles in two dimensions, the solution $u_N$ is obtained by collecting over all the local cell contributions $u_K$, representing the approximation of the solution $u$ on the cell $K$:
\begin{equation}
    u_K = u_\ell|_K.
\end{equation}
In order to be useful for local mesh refinement, the estimator $\eta_\ell$ is supposed to have a similar local structure, i.e., it is the sum of local contributions $\eta_K$ corresponding to $K$:
\begin{equation}\label{eq:eta_somme_des_etaK}
    \eta_\ell = \left( \sum_K \eta_K^{{2}} \right)^{{\frac{1}{2}}}. 
\end{equation}
The information encoded in the estimator map $\eta_K$ is used in different ways to {adapt or refine the mesh.} 
The main idea is to split locally a cell $K$, if the local error indicator $\eta_K$ is considered too large, see for instance \cite{pelle1996} or \cite{dorfler1996}. In a similar way one could use local coarsening, i.e., neighbouring cells are glued together to define a bigger cell where the local error indicator is small enough. The algorithm used to select cells for refinement (or coarsening) is called \textbf{Mark} 
and the step of mesh modification 
\textbf{Adapt}. 
%
%
The first natural question is about convergence of the adaptive algorithm (which does not follow from a priori error analysis). Under the condition of convergence, the second question is about the speed of convergence, measured in terms of number of unknowns.

Let us sketch a short argument for convergence. We consider the step from $\ell$ to $\ell+1$. For the sake of simplicity, we suppose that the boundary value problem \eqref{weak} corresponds to the minimisation of a quadratic energy functional. 
Moreover, we suppose that the finite element spaces 
$V_{\ell}$
form a sequence of conforming nested spaces ($V_{\ell} \subset V_{\ell + 1} \subset V$ ) 
and that there holds a Galerkin orthogonality
$(e_{\ell + 1},u_{\ell +1} -u_{\ell})_V = 0$, with $(\cdot,\cdot)_V$ the inner product in $V$. 
Then the errors $e_{\ell}:=\| u-u_{\ell}\|_V$ are related by
\[
e_{\ell+1}^2 =e_{\ell}^2 - \| u_{\ell+1} - u_{\ell}\|_V^2
\]
and we get geometrical convergence if for some constant $\beta>0$ 
\begin{equation}\label{eq:simplearg}
e_{\ell} \le \beta\| u_{\ell+1} - u_{\ell}\|_V.    
\end{equation}
since then $e_{\ell+1}\le \rho e_{\ell}$ with $\rho=\sqrt{1-1/\beta^2}$.

{
In order to establish (\ref{eq:simplearg}), it is convenient to assume the local lower bound 
\begin{equation}\label{eq:assumedlowerbound}
\eta_K \leq c \| u_{\ell+1} - u_{\ell} \|_{K}
\end{equation}
for all cells $K$ from the set of refined cells $\mathcal R_{\ell} \subset \mathcal K_{\ell}$ and to ensure that, for some $\theta>0$, the estimator linked to the cells in $\mathcal R_{\ell}$ is a least equal to a proportion $\theta$ of the total estimator:
\begin{equation}\label{eq:marking}
\sum_{K\in \mathcal R_{\ell}}\eta_K^2 \ge \theta^2 \sum_{K\in \mathcal K_{\ell}}\eta_K^2.
\end{equation}
Indeed, using \eqref{upperbound} and (\ref{eq:eta_somme_des_etaK}), then \eqref{eq:marking} and finally \eqref{eq:assumedlowerbound}, we get
$$
e_{\ell}^2 \leq C^2 \sum_{K\in \mathcal K_{\ell}}\eta_K^2 
           \leq \frac{C^2}{\theta^2} \sum_{K\in \mathcal R_{\ell}} \eta_K^2 
           \leq \frac{C^2 c^2}{\theta^2} \sum_{K\in \mathcal R_{\ell}} \| u_{\ell+1} - u_{\ell} \|_{K}^2 ,
$$
which implies \eqref{eq:simplearg} with $\beta = Cc/\theta.$
}


The assumed bound (\ref{eq:assumedlowerbound}) looks similar to what is called \emph{local efficiency}: there exists a constant $c>0$ independent of $N$, $u_N$ and $K$ such that, for every mesh cell $K$, there holds:
\begin{equation}
    \label{efficiencyloc}
    \eta_K \leq c \| u-u_N \|_{\omega_K},
\end{equation}
where $\omega_K$ is a collection of neighbouring cells of $K$ (patch) and 
$\| \cdot \|_{\omega_K}$ is an appropriate norm on $\omega_K$, see, e.g., 
\cite{verfurth2013}. This property strictly translates the fact that if the local discretization error is small, then the local estimator must also be small.
However, 
%
\eqref{efficiencyloc} does not imply~\eqref{eq:assumedlowerbound}
(the difficulty is not the presence of $\omega_K$ instead of $K$), and even (\ref{efficiencyloc}) might not always hold.
It might be thought that if an estimator does not satisfy this local efficiency property, mesh refinement could be misled. Indeed, regions where the discretization error is low could nevertheless be refined, as the local $\eta_K$ would be free to overestimate the actual discretization error. However, it turns out in practice that an estimator with 
poor efficiency can still be useful to drive the mesh refinement algorithm (see the discussion in \cite{carstensen2010}), and can even lead to optimal meshes. 

Condition (\ref{eq:marking}) states that a minimum percentage of the estimator contributions should give rise to refinement. It is usually called Dörfler marking or bulk chasing. Clearly, (\ref{eq:marking}) allows for a large number of possibilities for selecting~$\mathcal R_{\ell}$ and one naturally tends to refine the cells with the largest contributions. 
This choice is related to our second question about the speed of convergence. Choosing the cells to be refined as the largest contributors allows to bound the number of additional cells by the difference in estimators.
This idea basically leads to an estimate of the number of cells  needed to produce a given accuracy, see~\cite{stevenson2007} for the first proof.

The unavailability of lower bounds 
has led to the development of several more involved techniques, see section~\ref{sec:afem}.

\subsection{The goal-oriented paradigm}
\label{goal}
Many practical situations are focused on the approximation of a single quantity $J$, such as the drag- and lift-coefficients, average of heat transfer along a part of the boundary, stress intensity factors, a local stress or strain in elasticity, etc.
However, we observe that the aforementioned theory heavily relies on the functional analytic setting and the corresponding norms are used to measure the error $u-u_N$. Instead one is interested in
\[
| J(u)-J(u_N) |
\]
and the objective of goal-oriented error estimation is to provide estimators for this quantity. The main idea is to compute an approximation to the solution~$z$ to a \emph{dual} problem that involves $J$ in its right-hand side. This dual solution then allows to evaluate locally the sensitivity of $J$ with respect to the discretization error leading to 
an error estimator in the usual form, 
see for instance \cite{beckerrannacher2001},
and also 
\cite{becker1996,giles2002,gonzalez2014,han2005book,prudhomme1999,rognes2013}.







\section{Theory of adaptive finite element methods}\label{sec:afem} 

The mathematical theory of adaptive finite elements started with the work of Ivo Babu\v{s}ka in the 1970s (see \cite{babuska1978}), and found increasing interest with the works of Claes Johnson (see \cite{eriksson1985}), Mark Ainsworth and John Tinsley Oden (see \cite{ainsworth1993,ainsworth2000}), Pedro Morin, Ricardo Nochetto and Kunibert G. Siebert (see \cite{morin2000}), and Rüdiger Verfürth (see \cite{verfurth1999}) in the 1990s. Since then, it has made important theoretical progress, that we summarise below. 

\subsection{Optimality of adaptive algorithms} 

New insights have been brought to the theory of adaptive algorithms. 
{We first define what is meant by the notion of optimality.
 We assume a given 
procedure that, for any mesh $\mathcal K$, 
    computes an approximation $u_{\mathcal{K}}$ to the solution $u$ of a given partial differential equation, 
    usually by a finite element method based on approximation spaces $V_{\mathcal{K}}$ whose dimensions are denoted $N(\mathcal{K})$. We further assume possible 
    to approximate the solution $u$ by elements of $V_{\mathcal{K}}$ at a certain convergence rate $s>0$, which 
    mathematically translates into:
    there exists an absolute constant $C>0$ such that
    $$
    \inf_{v_{\mathcal{K}}\in V_{\mathcal{K}}, N(\mathcal{K})\le N}\left\| u-v_{\mathcal{K}}\right\| \le C N^{-s}
    \quad \forall N\in\mathbb N,
    $$
    in an appropriate norm. Then, an adaptive algorithm constructing a (in principle infinite) sequence of meshes $(\mathcal{K}_\ell)_{\ell=1,2,\cdots}$ and consequent solutions $(u_\ell)_{\ell=1,2,\cdots}$ is said to be
    (quasi-)optimal, if the sequence $(u_\ell)_{\ell=1,2,\cdots}$ has a similar error decay:
$$
    \left\| u-u_{\ell} \right\| \le C \left( N(\mathcal{K}_\ell) \right)^{-s}
    \quad \forall \ell \in\mathbb N,
$$
    (with possibly a different constant $C$)}.     
    One should mention that such an optimality result is much stronger than the classical \emph{a priori}  
    finite element theory, since no direct assumption on the regularity of the solution is made.
    

%

    This theory brings 
    the adaptive finite element technology
    close to the nonlinear approximation theory of multi-resolution wavelet methods, 
    see, e.g.,
    \cite{cohen2001} 
and \cite{karaivanov2003}. 
    It has been a very active research field, starting with the works of \cite{dorfler1996,morin2000}. The combination with nonlinear approximation theory has been achieved in \cite{binev2004,stevenson2007}, which lead to many further results for different methods/equations/algorithms, see for example \cite{cascon2008,BeckerMao08,Bonito10,Ferraz-LeiteOrtnerPraetorius10,Kreuzer18,Gantner22}.
   Now it has attained a state of maturity in the context of second-order linear elliptic partial differential equations, certain finite element methods and mesh refinement algorithms. The appearance of an underlying structure of the available results has lead to the Axiom-paper of \cite{carstensen2014}.
   Despite the large number of articles about application of goal-oriented adaptivity, theoretical results 
  {concerning their optimality} are still sparse. \cite{MommerStevenson09} provides the first optimality result, using the product of primal and dual estimators as an upper bound to the functional error; see also \cite{BeckerEstecahandyTrujillo11,FeischlPraetoriusZee16}.


\subsection{Other recent achievements}    
    
    Let us mention first that the described theory of 
    adaptive finite elements
    has been generalised to nonlinear monotone equations by \cite{FeischlFuhrerPraetorius14} and saddle-point systems as the Stokes equations, see \cite{BeckerMao11,bringmann2017bis,Feischl19,Feischl22}.
Also progress has been made in the understanding of parabolic problems, see \cite{Ern17}.
    
    Moreover, the theory has been generalised to control the overall work of the adaptive algorithm, and not just the dimension of the discrete spaces, see, e.g.,
    \cite{haberl2021,becker2023,Vohralik24}.
    This is particularly important from a practical viewpoint, since the various stopping criteria in the different nested loops can then be harmonised.
    

\section{Current engineering practice(s)} 

Now let us turn to the practical side of numerical methods, when they are used in practical computations in physics, applied sciences and industrial applications.
In this context, there have already been some initiatives to incorporate error control based on a posteriori error estimators as well as adaptive discretization techniques, see for instance the special volume
\cite{plewa2005} published earlier.
Nevertheless, the true potential of all these techniques seems to remain largely untapped.

\subsection{Accuracy vs. computational cost vs. complexity}

For numerical simulation oriented towards real applications, such as in the technological and industrial sector, the challenges are, at least, threefold:
\begin{enumerate}
    \item Simulation of complex systems, that may encompass, for instance, complex three-dimensional geometries, nonlinear partial differential equations with nonlinear boundary conditions, multi-physics or multi-scale phenomena, complex materials, etc.
    \item Simulation with reduced computation time and/or limited computational resources may be of major importance for certain applications, with the possible goal to even  provide predictions in real-time.
    \item Certified numerical simulation, with trustworthy and accurate predictions. 

\end{enumerate}
It is clear that the third point is related to error control and adaptive discretizations. In fact, the notions of verification and validation of finite element simulations have always been of importance within the numerical simulation community, see, e.g., 
\cite{babuska2004,babuska2006,babuska2007,stein2014,szabo2021}.
However this may be balanced or even be in conflict with the first two points, and this partly explains why error control is sometimes an undertaken issue, because it is not the top one priority for some practitioners devoted to a specific goal. 
{Our personal experience is that precision is often sacrificed for speed. One of our goals in these volumes is to show that these two requirements can be reconciled by the adaptivity and optimality provided by a posteriori error estimation.}

\subsection{Some initiatives}

Let us however mention a few examples in which a posteriori error estimations or control have been considered in an industrial context:
\begin{enumerate}
    \item For the neutron diffusion equation in nuclear reactor cores, see \cite{ciarletmadioneutron} (and see also \cite{taumhas2022impact} for an evaluation of the modelling error).
    \item For turbulent hydraulic and thermal-hydraulic simulations with a perspective of subsequent applications to the nuclear energy sector, see
    \cite{nassreddine2022,nassreddine2023} and \cite{dakroub2023}.
    \item Still with applications inspired by the nuclear industry, \cite{liu2017} address contact problems in elastostatics.    
    \item For real-time surgical simulation, see \cite{bui2017real}.
    as well as more details in section \ref{sub:realtimesurgery} below.
    \item Still for perspectives in surgical simulation and modelling of soft tissue, see \cite{duprez2020soft} and \cite{duprez:hal-04208610}.
    \item In the aeronautical sector, viscous compressible flows around airfoils were investigated in~\cite{Basileetal2022}
\end{enumerate}

Remark that these works are a step toward certified numerical simulation and have been highly motivated by a context where safety issues are fundamental.
Of course, while there are not many existing works, the above list is not exhaustive, and other illustrations are highlighted in some of the contributions to these volumes. 

\subsection{Some bottlenecks}

According to our own experience, we can point out the following bottlenecks relative to error control and adaptive discretization in routine industrial computations:

\begin{enumerate}
    \item For trustworthy certified numerical simulation, not only the discretization errors need to be controlled, but also modelling errors and uncertainties in parameter calibration. These latter can be of larger magnitude and it can be more involved to evaluate and reduce them.
    \item The design of a posteriori error estimators is, to a large extent, specific to the model under consideration, and it is challenging and time-consuming to design a posteriori error estimators well-suited for a complex model related to an industrial application.
    \item The efficient implementation of some error estimation techniques is not always straightforward in industrial or commercial codes.
    \item Error control and adaptive techniques are not always part of training related to numerical simulation, and a large part of the literature on the topic is not easily accessible to  outsiders.
\end{enumerate}

Note finally that part of the current activity in the field is inspired from these bottlenecks. Just to mention an example, various recent works have been dedicated to explore the interplay between adaptive finite element meshing and uncertainty quantification technologies, see \cite{bespalovI,bespalovII,eigelI,eigelII,guignardI,guignardII,oden2005}.

\section{Special topics and industrial applications}

Error estimation for numerical methods plays an important role in quality assurance. It enables users to better control the quality of their simulations. So-called non-standard numerical methods have also seen some advances in terms of error estimation. We will mention here a few numerical methods which offer alternatives to the standard finite element method based on Lagrange polynomials (see \ref{sub:galerkin}). We also discuss some of the most novel practical applications where a posteriori error estimates have been introduced, including surgical simulation and industrial scale fracture mechanics. 

\subsection{Real-time/interactive simulations}
\label{sub:realtimesurgery}

Because of the increasing need for interactive simulations, both for medical applications and robotic control and computer animations, but also in order to build digital twins of actual systems, the last 10 years or so have seen major progress towards real-time quality control of simulated solutions. 

This poses a number of problems which will be discussed in forthcoming special issues of Advances in Applied Mechanics. The first difficulty lies in computing the error estimate fast enough for the solution to be usable within the constraints of the practical application. For example, in surgical simulation, surgical guidance and surgical training, a range from 50 (for visual feedback) to 500 (for haptic feedback) solutions per second are necessary. Moreover, it is necessary to refine the mesh locally, based on the error estimate, within the time constraints posed by the application at hand. The interested reader can refer to the work of \cite{bui2018controlling,bui2017real,bui2019corotational}, which show the first real-time error estimation techniques for interactive mechanics simulations of bodies undergoing large deformations. 

Real-time error estimation for interactive deformable bodies still poses several challenges, especially as the demand for more realistic simulations in applications such as virtual surgery, computer animation, and virtual reality continues to grow. Here are some key challenges:

\begin{description}
    \item[Computational Efficiency] Developing error estimation techniques that can operate within the constraints of real-time simulation environments is crucial. These methods need to be computationally efficient to provide accurate estimates without significantly impacting the overall simulation performance.
\item[Integration with Simulation Frameworks] Integrating error estimation techniques seamlessly into existing simulation frameworks poses a challenge. These techniques need to be compatible with various simulation algorithms, such as finite element, meshfree, or position-based methods, to ensure wide\-spread adoption across different applications.
\item[Sensitivity to Simulation Parameters] Error estimation methods should be robust and reliable across different simulation scenarios and parameter settings. They should account for uncertainties and variations in material properties, boundary conditions, and external forces to provide accurate estimates in diverse environments. In the context of biomechanics, recent progress was made in a series of papers which attempt to disentangle model error from discretization error \cite{hauseux2017accelerating,hauseux2017calculating,hauseux2018quantifying}. 
In general, this connection between model and discretization error seems one of the richest direction of investigation, in particular when real-time data acquisition is possible and in the context of machine learning and artificial intelligence surrogate acceleration methods, e.g. \cite{deshpande2022probabilistic,deshpande2023convolution}.
\item[User Interaction and Control] Providing intuitive interfaces for users to interactively control error estimation and refinement processes is essential. Users should have the flexibility to adjust parameters, set error tolerances, and guide adaptive refinement strategies in real time to achieve desired simulation outcomes. A closer interaction between the human user and the computer could allow to switch between total computer control and total human control depending on the error level. 
\item[Guaranteed upper bounds] In real-time simulations of surgery, guaranteed upper bounds are crucial for ensuring the safety and effectiveness of the procedures. These bounds provide maximum limits on various parameters such as response time, computational complexity, and error rates. For instance, in robotic surgery, having guaranteed upper bounds on latency ensures that movements are executed within acceptable time frames to prevent accidents or errors. Similarly, computational complexity bounds ensure that simulations run efficiently without overwhelming hardware resources, allowing for smooth and accurate real-time feedback during surgical procedures. Overall, these guaranteed upper bounds are essential for maintaining reliability, accuracy, and safety in real-life applications like surgery simulations.
\end{description}

In short, one of the key remaining questions could be stated as: \emph{``Given multi-modal experimental observations obtained, in the best case, in real-time, and an underlying model, how sufficient is the data acquired to simulate the system ? How important and uncertain is the form of the model ? How important are the parameters used ? Are we better off simulating more scenarios (offline or online) or should we make more measurements ?''}. This direction seems ripe for fruitful investigations, as summarised in a recent review \cite{eftimie2023digital}.

\subsection{Strain smoothing}

Smoothed finite element methods see e.g. a review \cite{nguyen2008smooth} are based on the idea of transforming derivatives of shape functions into products with outward normals of smoothing cells. This suppresses the need to use Jacobian transformations, enables polytopal meshes, makes it possible to compute on extremely distorted meshes, alleviates numerical locking and can, in certain conditions provide upper bounds for the system’s energy. Several approaches were developed, mostly based on the methods shown in 
\cite{nguyen2008smooth}. 

Little work has been done in the area of error estimates for strain smoothing stabilized finite elements. The main question lies in the choice of the number of subcells. One subcell may provide equivalent methods to stress-based finite elements (dual), giving an upper bound to the energy. In the infinite limit, as the number of subcells goes to infinity, the standard displacement based (primal) finite element method is recovered. 
\cite{gonzalez2013efficient} presents ideas on how to obtain error indicators for smoothed finite element methods. 

\subsection{Partition of unity enrichment}

Extended finite element methods are partition of unity methods, similar to the generalized finite element method and a few others, including hp-clouds. The idea of the method 
\cite{bordas2023partition})
is to enrich the approximation by special functions which introduce known features of the exact solution. The idea is to improve the approximation property of the finite element (or meshfree) space by making it possible to reproduce the known features of the solution. See the  book \cite{bordas2023partition} and the recent review on error estimates for enriched approximations in \cite{gonzalez2023recovery}.

For example, to model a crack, a discontinuous function is introduced to take into account the jump across the crack surface. Asymptotic functions can be introduced to improve the accuracy of the solution close to the crack front. These methods were used in a wide variety of contexts, and called for specific treatments in terms of error estimation. Those are summarized in 
\cite{bordas2023partition}.  They were also implemented in a commercial code MorfeoCrack, based on the developments in the following seminal paper \cite{jin2017error}.

In short, standard recovery based methods fail because they are based on smoothing (projection on a polynomial space), whilst the enrichment functions are usually non-polynomial, see, e.g.,  \cite{bordas2007derivative}. This leads to smearing of the recovered solution, defeating its purpose to provide a higher quality solution to compare the raw finite element counterpart. Residual based methods have also to be adapted because of the extra terms present in the approximations, see
\cite{gerasimov2012,hild2009,ruter2013}. 

Several approaches exist for this, including subtracting the enrichment, see \cite{rodenas2008recovery}, or projecting on enriched spaces which are able to take into account the special features present in the solution space, see \cite{bordas2008simple,bordas2007derivative,duflot2008posteriori}. 
Other approaches were shown in the literature, such as \cite{panetier2010strict,prange2012error,ruter2013} and, recently, in \cite{bento2023recovery}, higher-order methods were investigated.

The goal-oriented techniques mentioned in 
\ref{goal} can be extended, and
a number of contributions exist in the field of goal-oriented error estimation 
for enriched approximations. The interested reader can refer to \cite{gonzalez2013efficient,gonzalez2014,gonzalez2021error,gonzalez2023recovery}.

The above papers introduce goal-oriented error estimates (GOEE), which quantify and control local errors in quantities of interest (QoI) for advanced engineering applications, such as aerospace, see also Section \ref{goal}. This includes a recovery-based error estimation technique for QoIs is presented, using an enhanced version of the Superconvergent Patch Recovery (SPR) technique, see \cite{zienkiewicz1987}. This approach provides nearly statically admissible stress fields, resulting in accurate estimations of local discretization error contributions to QoI. The technique requires reasonable computational cost and could be easily implemented into finite element codes or used as an independent postprocessing tool.

The error estimation in QoI relies on evaluating QoI through solving auxiliary  problems. Energy estimates are used to relate errors in QoI to the initial problem and auxiliary problem solutions. Explicit and implicit residual-based approaches and smoothing techniques for energy estimates can be used, as well as  the Zienkiewicz-Zhu estimate and SPR techniques. Bridging these approaches aims to obtain guaranteed upper bounds while retaining ease of implementation. The SPR-CX approach is   an efficient and simple goal-oriented adaptivity procedure for linear QoI in elasticity problems, extended to handle singular elasticity problems. 

\subsection{How important is it that the recovered (smoothed) solution satisfy the boundary conditions and the governing equations?}

We discuss here briefly the question of statical admissibility of recovered solutions, presented in the following: \cite{gonzalez2012role,bordas2023partition,gonzalez2023recovery}. These contributions  assesses the effect of statical admissibility and the ability of recovered solutions to represent singular solutions, along with the accuracy, local, and global effectivity of recovery-based error estimators for enriched finite element methods, such as the extended finite element method (XFEM). Two recovery techniques are studied: the superconvergent patch recovery procedure with equilibration and enrichment (SPR-CX) and the extended moving least squares recovery (XMLS). Both techniques enrich the recovered solutions, with SPR-CX enforcing equilibrium constraints. 

\emph{Numerical results highlight the necessity of extended recovery techniques in error estimators for this class of problems, with statically admissible recovered solutions yielding significant improvements in effectivities. This emphasizes the importance of both extended recovery procedures and statical admissibility for accurate assessment of the quality of enriched finite element approximations.}

\subsection{Meshfree/meshless methods}

Meshfree methods are numerical techniques used for solving partial differential equations without relying on a predefined mesh. Unlike traditional methods like finite element analysis, meshfree methods operate directly on scattered data points, making them particularly useful for problems with complex geometries or evolving domains. One notable advantage is their flexibility in handling irregular geometries and dynamic problems, which often pose challenges for mesh-based approaches.

A significant strength of meshfree methods lies in their ability to reproduce complex phenomena with high fidelity, especially in situations where traditional mesh-based methods struggle due to mesh distortion or excessive refinement requirements. Additionally, meshfree methods often exhibit good scalability and efficiency, particularly for problems involving large deformations or adaptive refinement. 

However, these methods also have their limitations. Reproducibility can be challenging since the results depend on the distribution of nodes or data points, which may vary between different simulations. Furthermore, achieving smooth solutions can be difficult, especially in regions with sparse data or irregular distributions of nodes.

For a comprehensive understanding of meshfree methods, one can refer to review papers such as \cite{nguyen2008meshless}, which provide an in-depth analysis of various aspects, including the theoretical foundations, implementation strategies, and applications of meshfree methods.

When it comes to error estimation in meshfree methods, several factors need to be taken into account. These include the choice of basis functions or shape functions, the accuracy of the numerical integration scheme, and the interpolation error associated with the scattered data points. Proper error estimation is crucial for assessing the reliability and accuracy of the computed solutions and guiding the refinement or adaptation strategies. Importantly, certain meshfree methods 
suffer from conditioning issues which may lead to algebraic errors, as discussed in this volume of Advances in Applied Mechanics, 
see \cite{papez2024algebraic}.
In particular this can apply to methods based on collocation, i.e., methods that, conversely to Galerkin techniques presented in \ref{sub:galerkin}, consist in writing the strong form directly at ``quasi-arbitrary'' sets of nodes within the domain, see \cite{jacquemin2023smart} for a review of adaptive schemes in meshless collocation, also known as the smart cloud method.

It is important to note that in meshfree methods, Galerkin orthogonality, a property commonly satisfied in traditional finite element methods, may not hold on the boundary. This is because the test functions employed in meshfree formulations typically do not vanish on the boundary, leading to deviations from orthogonality. Understanding and managing such deviations are essential for ensuring the accuracy and stability of the numerical solutions, particularly near domain boundaries.

Additionally, meshless methods do not employ polynomial approximations, which makes Gauss quadrature inexact. It is therefore critical to contain integration errors. This can be done by an interplay between using more background subcells in regions of interest as well as increasing the number of integration points per subcell, see \cite{meshless_integration_racz2012novel}.

The interested reader can refer to the following papers for additional references on error analysis for meshfree methods:  
\cite{ameshless_rnold1983asymptotic,
meshless_davydov2011adaptive,
meshless_dolbow1998introduction,
meshless_shells_adaptive_li2020adaptive,meshless_orkisz2008posteriori,meshless_park2003posteriori,
perazzo2008adaptive,
meshless_rabczuk2005adaptivity}, which 
includes applications to shell elements, optimal point placement in collocation methods, a posteriori error estimation driving point cloud adaptation, and applications to localized phenomena and large gradients. 

\subsection{Isogeometric analysis} 

Isogeometric analysis (IGA) is a computational technique that integrates the geometric design and analysis of structures or materials. Unlike traditional finite element methods, IGA employs the same basis functions to represent both the geometry and the solution field, typically using Non-Uniform Rational B-Splines (NURBS) or other spline-based representations. This seamless integration of geometric and simulation aspects offers several advantages over traditional methods.

One of the key strengths of IGA lies in its ability to precisely represent complex geometries using the same basis functions employed for the analysis, thereby eliminating the need for mesh generation and simplifying the workflow. This leads to significant reductions in pre-processing time and enables more accurate representations of geometric features.

Moreover, IGA facilitates the use of higher-order basis functions, which can lead to more accurate solutions, particularly for problems involving curved or irregular boundaries. By leveraging the smoothness properties of spline functions, IGA often produces smoother and more realistic solutions compared to traditional finite element methods.

However, despite these advantages, IGA also has its challenges. One notable limitation is the computational cost associated with the construction and manipulation of spline representations, especially for problems involving large-scale simulations or dynamic analyses. Additionally, the coupling between the geometry and analysis introduces complexities in the formulation and implementation of boundary conditions and geometric modifications, see, e.g., \cite{hu2018,nguyen2014}.

For a deeper understanding of IGA, interested individuals can refer to review papers such as \cite{nguyen2015isogeometric}, which provide comprehensive insights into the theoretical foundations, computational aspects, and applications of isogeometric analysis.

When it comes to error estimation in IGA, similar considerations apply as in traditional finite element methods, including the choice of basis functions, integration schemes, and interpolation errors. Proper error estimation is essential for assessing the accuracy and reliability of the computed solutions and guiding refinement strategies. 

Isogeometric analysis offers a powerful framework for integrating geometric design and analysis, with advantages in accuracy, efficiency, and geometric flexibility. However, challenges remain in terms of computational cost and the formulation of boundary conditions, highlighting the ongoing research efforts in this field.

The key difficulties in adaptive isogeometric simulations lie in the tensor-product nature of the underlying shape functions (NURBS), which, if nothing is done, leads to spurious propagation of refinements throughout the domain. Alternatives such as the geometry independent field approximation techniques avoid such issues by enabling the use of local refinement (through splines such as PHT splines, for example), whilst keeping the geometrical representation of the domain identical, and approximated by NURBS, see \cite{atroshchenko2018,IGA_GIFT_jansari2022adaptive,videla2019h,IGA_GIFT_videla2024shape}. 

Other exciting directions of research are provided in the following recent contributions including space-time analysis, as well as functional-type error estimates for isogeometric analysis, see \cite{langer2016posteriori,IGA_functional_matculevich2018functional,langer2019guaranteed}.





\section{Perspectives}

In \textit{Advances in Applied Mechanics (AAMS) Vol 58: Error control, adaptive discretizations, and applications, Part 1}, world-leading authors present cutting-edge research at the intersection of computational mechanics and applied mathematics, exploring innovative approaches to error control and adaptive discretizations across various fields.

\section*{Chapters}

\begin{enumerate}    
    \item \textbf{$hp$ adaptive Discontinuous Galerkin strategies driven by a posteriori error estimation with application to aeronautical flow problems} \citep{chapelier2024hp} presents recent developments about $h$, $p$, and $hp$-adaptive strategies driven by a posteriori error estimators using a high-order Discontinuous Galerkin finite element numerical framework, where $h$ stands for the mesh size and $p$ for the polynomial degree of the finite element approximation. A combination of error estimators tailored for the numerical method considered is presented along with smoothness indicators driving the decision to refine in $h$ or $p$. The $h$, $p$, and $hp$-adaptation strategies are described and applied to flow problems of aeronautical interest, including scale-resolving simulations of the transitional flow past a NACA0012 airfoil, scale-resolving simulations of the turbulent jet issued by a realistic nozzle geometry, and inviscid simulations of the transonic flow around a complex CRM aircraft geometry. For all cases considered, the interest of $h$, $p$, and $hp$-adaptation is demonstrated for easing the meshing process and increasing the resolution in flow regions of interest, enabling a significant reduction of the total number of computational degrees of freedom compared to manual meshing techniques and classical lower-order numerical approximations.

    \item \textbf{An anisotropic mesh adaptation method based on gradient recovery and optimal shape elements} \citep{fortin2023anisotropic} The authors present a complete mesh adaptation strategy applicable for controlling the discretization error on a finite element solution of Lagrange type of any degree. The method is described in detail, outlining the process of constructing a more accurate solution from a finite element solution using a gradient recovery method. The error is then estimated as the difference between the reconstructed and the initial solution. Subsequently, the mesh is modified using local operations to minimize the error on the gradient. Additionally, a number of numerical examples are provided to illustrate the effectiveness of the approach.

    \item \textbf{Model reduction techniques for parametrized nonlinear partial differential equations} \citep{nguyen2024model} The authors present model reduction techniques for parametrized nonlinear partial differential equations (PDEs). The main ingredients of their approach include reduced basis (RB) spaces to provide rapidly convergent approximations to the parametric manifold, Galerkin projection of the underlying PDEs onto the RB space to reduce the number of degrees of freedom, and empirical interpolation schemes to rapidly evaluate the nonlinear terms associated with the Galerkin projection. They devise a first-order empirical interpolation method to construct an inexpensive and stable interpolation of the nonlinear terms. Two different hyper-reduction strategies are considered: hyper-reduction followed by linearization, and linearization followed by hyper-reduction. The authors extend empirical interpolation to nonintrusive model reduction and apply it to compressible flows in both supersonic and hypersonic regimes. They present numerical results to illustrate the accuracy, efficiency, and stability of the reduced-order models. The interested reader can refer to \cite{kerfriden2014certification,hoang2016fast,goury2016automatised,agathos2020parametrized,hoang2021domain,hoang2022projection,chen2024reduced} for adaptive (certified) model order reduction for non-linear problems and localised phenomena such as fracture, molecular dynamics and inverse problems.
    
    \item \textbf{Adaptive finite elements for obstacle problems} \citep{gustafsson2024adaptive} The author summarize three applications of the obstacle problem to membrane contact, elastoplastic torsion, and cavitation modelling, demonstrating how the resulting models can be solved using mixed finite elements. He highlights the challenge of constructing fixed computational meshes for any inequality-constrained problem due to the unknown shape of the coincidence set. Consequently, he demonstrates how $h$--adaptivity can be utilized to resolve the unknown coincidence set. Additionally, the author discusses practical challenges that must be overcome in the application of the adaptive methods.

    \item \textbf{A Posteriori Error Identities and Estimates of modelling Errors} \citep{repin2024posteriori}  This Chapter discusses a posteriori estimation methods for mathematical models based on partial differential equations. The analysis is based on functional identities of a special kind, reflecting the most general relations that hold for deviations from the exact solution of a boundary value problem. These identities do not depend on special properties of approximations and contain no mesh-dependent constants, making them valid for any function in the admissible (energy) class. This universality enables the control of accuracy for various numerical approximations and the comparison of solutions of mathematical models. The capabilities are demonstrated with the paradigm of modelling errors generated by simplifications of the original problem. Three groups of problems are discussed, where errors of simplification have different origins: errors arising from simplifying coefficients of the equation, errors associated with simplifying geometry, and dimension reduction errors. It is shown that in any of these cases, the desired error estimates follow from the general a posteriori identity after proper specification of the functional spaces and operators associated with the boundary value problem.

    \item \textbf{Exact error control for variational problems via convex duality and explicit flux reconstruction} \citep{bartels2024exact}
    A posteriori error estimates serve as crucial tools for bounding discretization errors in terms of computable quantities, bypassing the need for establishing often challenging regularity conditions. However, for problems involving non-linearity, non-differentiability, jumping coefficients, or finite element methods with anisotropic triangulations, such estimates can often result in large factors, leading to sub-optimal error estimates. To address this issue, exact and explicit error representations are derived using convex duality arguments, effectively avoiding such effects.

    \item \textbf{Algebraic error in numerical PDEs and its estimation} \citep{papez2024algebraic}  The paper discusses error estimation for the Poisson model problem, focusing on residual-based and flux reconstruction error estimates. Residual-based estimates, while useful for estimating discretization errors, have limitations due to the lack of exact algebraic solutions. The paper proposes adjustments to these estimates but highlights drawbacks. Flux reconstruction offers local error indicators and guaranteed bounds but is computationally intensive. The paper emphasizes the importance of reducing algebraic error for accurate discretization error estimation and suggests further research on adaptive mesh refinement and algebraic solver stopping criteria. It also advocates for rigorous methods to mitigate algebraic errors in various numerical methods. It is worth noting that other methods including meshfree methods, enriched partition of unity methods, see \cite{agathos2016well,agathos2016stable,agathos2018stable,agathos2019unified,agathos2019improving,agathos2020parametrized,bordas2023partition}, isogeometric analysis, see \cite{langer2016posteriori,yu2018adaptive,IGA_functional_matculevich2018functional,videla2019h,langer2019guaranteed,jansari2022adaptive} and collocation methods, see \cite{ameshless_rnold1983asymptotic,jia2019adaptive}, smart point clouds, see \cite{perazzo2008adaptive,jacquemin2023smart}, are also concerned by algebraic errors due to ill-conditioning of the stiffness matrix. The interested reader can refer to the following papers. 

\end{enumerate}

\bibliographystyle{elsarticle-harv}
\bibliography{biblio}

\begin{thebibliography}{182}
\expandafter\ifx\csname natexlab\endcsname\relax\def\natexlab#1{#1}\fi
\providecommand{\url}[1]{\texttt{#1}}
\providecommand{\href}[2]{#2}
\providecommand{\path}[1]{#1}
\providecommand{\DOIprefix}{doi:}
\providecommand{\ArXivprefix}{arXiv:}
\providecommand{\URLprefix}{URL: }
\providecommand{\Pubmedprefix}{pmid:}
\providecommand{\doi}[1]{\href{http://dx.doi.org/#1}{\path{#1}}}
\providecommand{\Pubmed}[1]{\href{pmid:#1}{\path{#1}}}
\providecommand{\bibinfo}[2]{#2}
\ifx\xfnm\relax \def\xfnm[#1]{\unskip,\space#1}\fi
\bibitem[{Agathos et~al.(2019a)Agathos, Bordas and
  Chatzi}]{agathos2019improving}
\bibinfo{author}{Agathos, K.}, \bibinfo{author}{Bordas, S.P.},
  \bibinfo{author}{Chatzi, E.}, \bibinfo{year}{2019}a.
\newblock \bibinfo{title}{Improving the conditioning of xfem/gfem for fracture
  mechanics problems through enrichment quasi-orthogonalization}.
\newblock \bibinfo{journal}{Computer Methods in Applied Mechanics and
  Engineering} \bibinfo{volume}{346}, \bibinfo{pages}{1051--1073}.
\bibitem[{Agathos et~al.(2020)Agathos, Bordas and
  Chatzi}]{agathos2020parametrized}
\bibinfo{author}{Agathos, K.}, \bibinfo{author}{Bordas, S.P.},
  \bibinfo{author}{Chatzi, E.}, \bibinfo{year}{2020}.
\newblock \bibinfo{title}{Parametrized reduced order modeling for cracked
  solids}.
\newblock \bibinfo{journal}{International Journal for Numerical Methods in
  Engineering} \bibinfo{volume}{121}, \bibinfo{pages}{4537--4565}.
\bibitem[{Agathos et~al.(2016a)Agathos, Chatzi and Bordas}]{agathos2016stable}
\bibinfo{author}{Agathos, K.}, \bibinfo{author}{Chatzi, E.},
  \bibinfo{author}{Bordas, S.P.}, \bibinfo{year}{2016}a.
\newblock \bibinfo{title}{Stable 3d extended finite elements with higher order
  enrichment for accurate non planar fracture}.
\newblock \bibinfo{journal}{Computer Methods in Applied Mechanics and
  Engineering} \bibinfo{volume}{306}, \bibinfo{pages}{19--46}.
\bibitem[{Agathos et~al.(2019b)Agathos, Chatzi and Bordas}]{agathos2019unified}
\bibinfo{author}{Agathos, K.}, \bibinfo{author}{Chatzi, E.},
  \bibinfo{author}{Bordas, S.P.}, \bibinfo{year}{2019}b.
\newblock \bibinfo{title}{A unified enrichment approach addressing blending and
  conditioning issues in enriched finite elements}.
\newblock \bibinfo{journal}{Computer Methods in Applied Mechanics and
  Engineering} \bibinfo{volume}{349}, \bibinfo{pages}{673--700}.
\bibitem[{Agathos et~al.(2016b)Agathos, Chatzi, Bordas and
  Talaslidis}]{agathos2016well}
\bibinfo{author}{Agathos, K.}, \bibinfo{author}{Chatzi, E.},
  \bibinfo{author}{Bordas, S.P.}, \bibinfo{author}{Talaslidis, D.},
  \bibinfo{year}{2016}b.
\newblock \bibinfo{title}{A well-conditioned and optimally convergent xfem for
  3d linear elastic fracture}.
\newblock \bibinfo{journal}{International Journal for Numerical Methods in
  Engineering} \bibinfo{volume}{105}, \bibinfo{pages}{643--677}.
\bibitem[{Agathos et~al.(2018)Agathos, Ventura, Chatzi and
  Bordas}]{agathos2018stable}
\bibinfo{author}{Agathos, K.}, \bibinfo{author}{Ventura, G.},
  \bibinfo{author}{Chatzi, E.}, \bibinfo{author}{Bordas, S.P.},
  \bibinfo{year}{2018}.
\newblock \bibinfo{title}{Stable 3d xfem/vector level sets for non-planar 3d
  crack propagation and comparison of enrichment schemes}.
\newblock \bibinfo{journal}{International Journal for Numerical Methods in
  Engineering} \bibinfo{volume}{113}, \bibinfo{pages}{252--276}.
\bibitem[{Ainsworth and Oden(1993)}]{ainsworth1993}
\bibinfo{author}{Ainsworth, M.}, \bibinfo{author}{Oden, J.T.},
  \bibinfo{year}{1993}.
\newblock \bibinfo{title}{A unified approach to a posteriori error estimation
  using element residual methods}.
\newblock \bibinfo{journal}{Numerische Mathematik} \bibinfo{volume}{65},
  \bibinfo{pages}{23--50}.
\bibitem[{Ainsworth and Oden(2000)}]{ainsworth2000}
\bibinfo{author}{Ainsworth, M.}, \bibinfo{author}{Oden, J.T.},
  \bibinfo{year}{2000}.
\newblock \bibinfo{title}{A posteriori error estimation in finite element
  analysis}.
\newblock Pure Appl. Math., Wiley-Intersci. Ser. Texts Monogr. Tracts,
  \bibinfo{publisher}{Chichester: Wiley}.
\bibitem[{Arnold(1982)}]{arnold-82}
\bibinfo{author}{Arnold, D.N.}, \bibinfo{year}{1982}.
\newblock \bibinfo{title}{An interior penalty finite element method with
  discontinuous elements}.
\newblock \bibinfo{journal}{SIAM Journal on Numerical Analysis}
  \bibinfo{volume}{19}, \bibinfo{pages}{742--760}.
\bibitem[{Arnold and Wendland(1983)}]{ameshless_rnold1983asymptotic}
\bibinfo{author}{Arnold, D.N.}, \bibinfo{author}{Wendland, W.L.},
  \bibinfo{year}{1983}.
\newblock \bibinfo{title}{On the asymptotic convergence of collocation
  methods}.
\newblock \bibinfo{journal}{Mathematics of Computation} \bibinfo{volume}{41},
  \bibinfo{pages}{349--381}.
\bibitem[{Atroshchenko et~al.(2018)Atroshchenko, Tomar, Xu and
  Bordas}]{atroshchenko2018}
\bibinfo{author}{Atroshchenko, E.}, \bibinfo{author}{Tomar, S.},
  \bibinfo{author}{Xu, G.}, \bibinfo{author}{Bordas, S.P.A.},
  \bibinfo{year}{2018}.
\newblock \bibinfo{title}{Weakening the tight coupling between geometry and
  simulation in isogeometric analysis: from sub- and super-geometric analysis
  to geometry-independent field approximation ({GIFT})}.
\newblock \bibinfo{journal}{International Journal for Numerical Methods in
  Engineering} \bibinfo{volume}{114}, \bibinfo{pages}{1131--1159}.
\bibitem[{Babu{\v{s}}ka et~al.(2007)Babu{\v{s}}ka, Nobile and
  Tempone}]{babuska2007}
\bibinfo{author}{Babu{\v{s}}ka, I.}, \bibinfo{author}{Nobile, F.},
  \bibinfo{author}{Tempone, R.}, \bibinfo{year}{2007}.
\newblock \bibinfo{title}{Reliability of computational science}.
\newblock \bibinfo{journal}{Numerical Methods for Partial Differential
  Equations} \bibinfo{volume}{23}, \bibinfo{pages}{753--784}.
\bibitem[{Babuska and Oden(2004)}]{babuska2004}
\bibinfo{author}{Babuska, I.}, \bibinfo{author}{Oden, J.T.},
  \bibinfo{year}{2004}.
\newblock \bibinfo{title}{Verification and validation in computational
  engineering and science: basic concepts}.
\newblock \bibinfo{journal}{Computer Methods in Applied Mechanics and
  Engineering} \bibinfo{volume}{193}, \bibinfo{pages}{4057--4066}.
\bibitem[{Babu{\v{s}}ka and Oden(2006)}]{babuska2006}
\bibinfo{author}{Babu{\v{s}}ka, I.}, \bibinfo{author}{Oden, J.T.},
  \bibinfo{year}{2006}.
\newblock \bibinfo{title}{The reliability of computer predictions: can they be
  trusted?}
\newblock \bibinfo{journal}{International Journal of Numerical Analysis and
  Modeling} \bibinfo{volume}{3}, \bibinfo{pages}{255--272}.
\bibitem[{Babu{\v{s}}ka and Rheinboldt(1978)}]{babuska1978}
\bibinfo{author}{Babu{\v{s}}ka, I.}, \bibinfo{author}{Rheinboldt, W.C.},
  \bibinfo{year}{1978}.
\newblock \bibinfo{title}{Error estimates for adaptive finite element
  computations}.
\newblock \bibinfo{journal}{SIAM Journal on Numerical Analysis}
  \bibinfo{volume}{15}, \bibinfo{pages}{736--754}.
\bibitem[{Babu\v{s}ka et~al.(1992)Babu\v{s}ka, Dur\'{a}n and
  Rodr\'{\i}guez}]{babuska1992}
\bibinfo{author}{Babu\v{s}ka, I.}, \bibinfo{author}{Dur\'{a}n, R.},
  \bibinfo{author}{Rodr\'{\i}guez, R.}, \bibinfo{year}{1992}.
\newblock \bibinfo{title}{Analysis of the efficiency of an a posteriori error
  estimator for linear triangular finite elements}.
\newblock \bibinfo{journal}{SIAM Journal on Numerical Analysis}
  \bibinfo{volume}{29}, \bibinfo{pages}{947--964}.
\bibitem[{Bank and Weiser(1985)}]{bank1985}
\bibinfo{author}{Bank, R.E.}, \bibinfo{author}{Weiser, A.},
  \bibinfo{year}{1985}.
\newblock \bibinfo{title}{Some a posteriori error estimators for elliptic
  partial differential equations}.
\newblock \bibinfo{journal}{Mathematics of Computation} \bibinfo{volume}{44},
  \bibinfo{pages}{283--301}.
\bibitem[{Bartels and Kaltenbach(2024)}]{bartels2024exact}
\bibinfo{author}{Bartels, S.}, \bibinfo{author}{Kaltenbach, A.},
  \bibinfo{year}{2024}.
\newblock \bibinfo{title}{Exact error control for variational problems via
  convex duality and explicit flux reconstruction}.
\newblock \bibinfo{journal}{Advances in Applied Mechanics (AAMS)}
  \bibinfo{volume}{58}.
\bibitem[{Basile et~al.(2022)Basile, Chapelier, de~la Llave~Plata, Laraufie and
  Frey}]{Basileetal2022}
\bibinfo{author}{Basile, F.}, \bibinfo{author}{Chapelier, J.B.},
  \bibinfo{author}{de~la Llave~Plata, M.}, \bibinfo{author}{Laraufie, R.},
  \bibinfo{author}{Frey, P.}, \bibinfo{year}{2022}.
\newblock \bibinfo{title}{Unstructured {{\(h\)}}- and \emph{hp}-adaptive
  strategies for discontinuous galerkin methods based on \emph{a posteriori}
  error estimation for compressible flows}.
\newblock \bibinfo{journal}{Comput. Fluids} \bibinfo{volume}{233},
  \bibinfo{pages}{21}.
\newblock \DOIprefix\doi{10.1016/j.compfluid.2021.105245}. \bibinfo{note}{id/No
  105245}.
\bibitem[{Becker et~al.(2023)Becker, Brunner, Innerberger, Melenk and
  Praetorius}]{becker2023}
\bibinfo{author}{Becker, R.}, \bibinfo{author}{Brunner, M.},
  \bibinfo{author}{Innerberger, M.}, \bibinfo{author}{Melenk, J.M.},
  \bibinfo{author}{Praetorius, D.}, \bibinfo{year}{2023}.
\newblock \bibinfo{title}{Cost-optimal adaptive iterative linearized {FEM} for
  semilinear elliptic {PDEs}}.
\newblock \bibinfo{journal}{European Series in Applied and Industrial
  Mathematics (ESAIM): Mathematical Modelling and Numerical Analysis}
  \bibinfo{volume}{57}, \bibinfo{pages}{2193--2225}.
\bibitem[{Becker et~al.(2011)Becker, Estecahandy and
  Trujillo}]{BeckerEstecahandyTrujillo11}
\bibinfo{author}{Becker, R.}, \bibinfo{author}{Estecahandy, E.},
  \bibinfo{author}{Trujillo, D.}, \bibinfo{year}{2011}.
\newblock \bibinfo{title}{Weighted marking for goal-oriented adaptive finite
  element methods}.
\newblock \bibinfo{journal}{SIAM Journal on Numerical Analysis}
  \bibinfo{volume}{49}, \bibinfo{pages}{2451--2469}.
\bibitem[{Becker and Mao(2008)}]{BeckerMao08}
\bibinfo{author}{Becker, R.}, \bibinfo{author}{Mao, S.}, \bibinfo{year}{2008}.
\newblock \bibinfo{title}{An optimally convergent adaptive mixed finite element
  method}.
\newblock \bibinfo{journal}{Numerische Mathematik} \bibinfo{volume}{111},
  \bibinfo{pages}{35--54}.
\bibitem[{Becker and Mao(2011)}]{BeckerMao11}
\bibinfo{author}{Becker, R.}, \bibinfo{author}{Mao, S.}, \bibinfo{year}{2011}.
\newblock \bibinfo{title}{Quasi-optimality of adaptive non-conforming finite
  element methods for the {S}tokes equations}.
\newblock \bibinfo{journal}{SIAM Journal on Numerical Analysis}
  \bibinfo{volume}{49}, \bibinfo{pages}{970--991}.
\bibitem[{Becker and Rannacher(1996)}]{becker1996}
\bibinfo{author}{Becker, R.}, \bibinfo{author}{Rannacher, R.},
  \bibinfo{year}{1996}.
\newblock \bibinfo{title}{A feed-back approach to error control in finite
  element methods: {Basic} analysis and examples}.
\newblock \bibinfo{journal}{East-West Journal of Numerical Mathematics}
  \bibinfo{volume}{4}, \bibinfo{pages}{237--264}.
\bibitem[{Becker and Rannacher(2001)}]{beckerrannacher2001}
\bibinfo{author}{Becker, R.}, \bibinfo{author}{Rannacher, R.},
  \bibinfo{year}{2001}.
\newblock \bibinfo{title}{An optimal control approach to a posteriori error
  estimation in finite element methods}.
\newblock \bibinfo{journal}{Acta Numerica} \bibinfo{volume}{10},
  \bibinfo{pages}{1--102}.
\bibitem[{{Beir\~{a}o~da~Veiga} et~al.(2013){Beir\~{a}o~da~Veiga}, Brezzi,
  Cangiani, Manzini, Marini and Russo}]{beirao2013vem}
\bibinfo{author}{{Beir\~{a}o~da~Veiga}, L.}, \bibinfo{author}{Brezzi, F.},
  \bibinfo{author}{Cangiani, A.}, \bibinfo{author}{Manzini, G.},
  \bibinfo{author}{Marini, L.D.}, \bibinfo{author}{Russo, A.},
  \bibinfo{year}{2013}.
\newblock \bibinfo{title}{Basic principles of virtual element methods}.
\newblock \bibinfo{journal}{Mathematical Models and Methods in Applied
  Sciences} \bibinfo{volume}{23}, \bibinfo{pages}{199--214}.
\bibitem[{Bento et~al.(2023)Bento, Proen{\c{c}}a and
  Duarte}]{bento2023recovery}
\bibinfo{author}{Bento, M.H.}, \bibinfo{author}{Proen{\c{c}}a, S.P.},
  \bibinfo{author}{Duarte, C.A.}, \bibinfo{year}{2023}.
\newblock \bibinfo{title}{Recovery strategies, a posteriori error estimation,
  and local error indication for second-order g/xfem and fem}.
\newblock \bibinfo{journal}{International Journal for Numerical Methods in
  Engineering} \bibinfo{volume}{124}, \bibinfo{pages}{3025--3062}.
\bibitem[{Bernardi and Maday(1997)}]{bernardi1997}
\bibinfo{author}{Bernardi, C.}, \bibinfo{author}{Maday, Y.},
  \bibinfo{year}{1997}.
\newblock \bibinfo{title}{Spectral methods}, in: \bibinfo{booktitle}{Handbook
  of numerical analysis, {V}ol. {V}}. \bibinfo{publisher}{North-Holland,
  Amsterdam}. volume~\bibinfo{volume}{V} of \textit{\bibinfo{series}{Handb.
  Numer. Anal.}}, pp. \bibinfo{pages}{209--485}.
\bibitem[{Bertoluzza(1995)}]{bertoluzza1995}
\bibinfo{author}{Bertoluzza, S.}, \bibinfo{year}{1995}.
\newblock \bibinfo{title}{A posteriori error estimates for the wavelet
  {Galerkin} method}.
\newblock \bibinfo{journal}{Applied Mathematics Letters} \bibinfo{volume}{8},
  \bibinfo{pages}{1--6}.
\bibitem[{Bertoluzza et~al.(1994)Bertoluzza, Naldi and Ravel}]{bertoluzza1994}
\bibinfo{author}{Bertoluzza, S.}, \bibinfo{author}{Naldi, G.},
  \bibinfo{author}{Ravel, J.C.}, \bibinfo{year}{1994}.
\newblock \bibinfo{title}{Wavelet methods for the numerical solution of
  boundary value problems on the interval}, in: \bibinfo{booktitle}{Wavelets:
  theory, algorithms, and applications. Proceedings of the international
  conference on wavelets, held in Taormina, Italy, October 14-20, 1993}.
  \bibinfo{publisher}{San Diego, CA: Academic Press}, pp.
  \bibinfo{pages}{425--448}.
\bibitem[{Bespalov and Silvester(2023)}]{bespalovII}
\bibinfo{author}{Bespalov, A.}, \bibinfo{author}{Silvester, D.},
  \bibinfo{year}{2023}.
\newblock \bibinfo{title}{Error estimation and adaptivity for stochastic
  collocation finite elements {P}art {II}: {M}ultilevel approximation}.
\newblock \bibinfo{journal}{SIAM Journal on Scientific Computing}
  \bibinfo{volume}{45}, \bibinfo{pages}{A784--A800}.
\bibitem[{Bespalov et~al.(2022)Bespalov, Silvester and Xu}]{bespalovI}
\bibinfo{author}{Bespalov, A.}, \bibinfo{author}{Silvester, D.J.},
  \bibinfo{author}{Xu, F.}, \bibinfo{year}{2022}.
\newblock \bibinfo{title}{Error estimation and adaptivity for stochastic
  collocation finite elements {P}art {I}: {S}ingle-level approximation}.
\newblock \bibinfo{journal}{SIAM Journal on Scientific Computing}
  \bibinfo{volume}{44}, \bibinfo{pages}{A3393--A3412}.
\bibitem[{Binev et~al.(2004)Binev, Dahmen and DeVore}]{binev2004}
\bibinfo{author}{Binev, P.}, \bibinfo{author}{Dahmen, W.},
  \bibinfo{author}{DeVore, R.}, \bibinfo{year}{2004}.
\newblock \bibinfo{title}{Adaptive finite element methods with convergence
  rates}.
\newblock \bibinfo{journal}{Numerische Mathematik} \bibinfo{volume}{97},
  \bibinfo{pages}{219--268}.
\bibitem[{Bonito and Nochetto(2010)}]{Bonito10}
\bibinfo{author}{Bonito, A.}, \bibinfo{author}{Nochetto, R.H.},
  \bibinfo{year}{2010}.
\newblock \bibinfo{title}{Quasi-optimal convergence rate of an adaptive
  discontinuous {G}alerkin method}.
\newblock \bibinfo{journal}{SIAM Journal on Numerical Analysis}
  \bibinfo{volume}{48}, \bibinfo{pages}{734--771}.
\bibitem[{Bordas and Duflot(2007)}]{bordas2007derivative}
\bibinfo{author}{Bordas, S.}, \bibinfo{author}{Duflot, M.},
  \bibinfo{year}{2007}.
\newblock \bibinfo{title}{Derivative recovery and a posteriori error estimate
  for extended finite elements}.
\newblock \bibinfo{journal}{Computer Methods in Applied Mechanics and
  Engineering} \bibinfo{volume}{196}, \bibinfo{pages}{3381--3399}.
\bibitem[{Bordas et~al.(2008)Bordas, Duflot and Le}]{bordas2008simple}
\bibinfo{author}{Bordas, S.}, \bibinfo{author}{Duflot, M.},
  \bibinfo{author}{Le, P.}, \bibinfo{year}{2008}.
\newblock \bibinfo{title}{A simple error estimator for extended finite
  elements}.
\newblock \bibinfo{journal}{Communications in Numerical Methods in Engineering}
  \bibinfo{volume}{24}, \bibinfo{pages}{961--971}.
\bibitem[{Bordas and Menk(2023)}]{bordas2023partition}
\bibinfo{author}{Bordas, S.}, \bibinfo{author}{Menk, A.}, \bibinfo{year}{2023}.
\newblock \bibinfo{title}{Partition of Unity Methods}.
\newblock \bibinfo{publisher}{Wiley}.
\bibitem[{Braess and Sch\"{o}berl(2008)}]{braess2008}
\bibinfo{author}{Braess, D.}, \bibinfo{author}{Sch\"{o}berl, J.},
  \bibinfo{year}{2008}.
\newblock \bibinfo{title}{Equilibrated residual error estimator for edge
  elements}.
\newblock \bibinfo{journal}{Mathematics of Computation} \bibinfo{volume}{77},
  \bibinfo{pages}{651--672}.
\bibitem[{Brenner and Scott(2008)}]{brenner2008}
\bibinfo{author}{Brenner, S.C.}, \bibinfo{author}{Scott, L.R.},
  \bibinfo{year}{2008}.
\newblock \bibinfo{title}{The mathematical theory of finite element methods}.
  volume~\bibinfo{volume}{15} of \textit{\bibinfo{series}{Texts in Applied
  Mathematics}}.
\newblock \bibinfo{edition}{Third} ed., \bibinfo{publisher}{Springer, New
  York}.
\bibitem[{Bringmann and Carstensen(2017)}]{bringmann2017bis}
\bibinfo{author}{Bringmann, P.}, \bibinfo{author}{Carstensen, C.},
  \bibinfo{year}{2017}.
\newblock \bibinfo{title}{An adaptive least-squares {FEM} for the {Stokes}
  equations with optimal convergence rates}.
\newblock \bibinfo{journal}{Numerische Mathematik} \bibinfo{volume}{135},
  \bibinfo{pages}{459--492}.
\bibitem[{Buffa et~al.(2022)Buffa, Gantner, Giannelli, Praetorius and
  V\'{a}zquez}]{Gantner22}
\bibinfo{author}{Buffa, A.}, \bibinfo{author}{Gantner, G.},
  \bibinfo{author}{Giannelli, C.}, \bibinfo{author}{Praetorius, D.},
  \bibinfo{author}{V\'{a}zquez, R.}, \bibinfo{year}{2022}.
\newblock \bibinfo{title}{Mathematical foundations of adaptive isogeometric
  analysis}.
\newblock \bibinfo{journal}{Archives of Computational Methods in Engineering}
  \bibinfo{volume}{29}, \bibinfo{pages}{4479--4555}.
\bibitem[{Bui et~al.(2023)Bui, Duprez, Rohan, Lejeune, Bordas, Bucki and
  Chouly}]{duprez:hal-04208610}
\bibinfo{author}{Bui, H.P.}, \bibinfo{author}{Duprez, M.},
  \bibinfo{author}{Rohan, P.Y.}, \bibinfo{author}{Lejeune, A.},
  \bibinfo{author}{Bordas, S.P.A.}, \bibinfo{author}{Bucki, M.},
  \bibinfo{author}{Chouly, F.}, \bibinfo{year}{2023}.
\newblock \bibinfo{title}{{Automatic mesh refinement for soft tissue}}.
\newblock \URLprefix \url{https://hal.science/hal-04208610}. \bibinfo{note}{hAL
  preprint}.
\bibitem[{Bui et~al.(2019)Bui, Tomar and Bordas}]{bui2019corotational}
\bibinfo{author}{Bui, H.P.}, \bibinfo{author}{Tomar, S.},
  \bibinfo{author}{Bordas, S.P.}, \bibinfo{year}{2019}.
\newblock \bibinfo{title}{Corotational cut finite element method for real-time
  surgical simulation: Application to needle insertion simulation}.
\newblock \bibinfo{journal}{Computer Methods in Applied Mechanics and
  Engineering} \bibinfo{volume}{345}, \bibinfo{pages}{183--211}.
\bibitem[{Bui et~al.(2018a)Bui, Tomar, Courtecuisse, Audette, Cotin and
  Bordas}]{bui2018controlling}
\bibinfo{author}{Bui, H.P.}, \bibinfo{author}{Tomar, S.},
  \bibinfo{author}{Courtecuisse, H.}, \bibinfo{author}{Audette, M.},
  \bibinfo{author}{Cotin, S.}, \bibinfo{author}{Bordas, S.P.},
  \bibinfo{year}{2018}a.
\newblock \bibinfo{title}{Controlling the error on target motion through
  real-time mesh adaptation: applications to deep brain stimulation}.
\newblock \bibinfo{journal}{International Journal for Numerical Methods in
  Biomedical Engineering} \bibinfo{volume}{34}, \bibinfo{pages}{e2958}.
\bibitem[{Bui et~al.(2018b)Bui, Tomar, Courtecuisse, Cotin and
  Bordas}]{bui2017real}
\bibinfo{author}{Bui, H.P.}, \bibinfo{author}{Tomar, S.},
  \bibinfo{author}{Courtecuisse, H.}, \bibinfo{author}{Cotin, S.},
  \bibinfo{author}{Bordas, S.P.A.}, \bibinfo{year}{2018}b.
\newblock \bibinfo{title}{Real-time error control for surgical simulation}.
\newblock \bibinfo{journal}{IEEE Transactions on Biomedical Engineering}
  \bibinfo{volume}{65}, \bibinfo{pages}{596--607}.
\bibitem[{Burman et~al.(2015)Burman, Claus, Hansbo, Larson and
  Massing}]{burman2015}
\bibinfo{author}{Burman, E.}, \bibinfo{author}{Claus, S.},
  \bibinfo{author}{Hansbo, P.}, \bibinfo{author}{Larson, M.G.},
  \bibinfo{author}{Massing, A.}, \bibinfo{year}{2015}.
\newblock \bibinfo{title}{Cut{FEM}: discretizing geometry and partial
  differential equations}.
\newblock \bibinfo{journal}{International Journal for Numerical Methods in
  Engineering} \bibinfo{volume}{104}, \bibinfo{pages}{472--501}.
\bibitem[{Canuto et~al.(2006)Canuto, Hussaini, Quarteroni and
  Zang}]{canuto2006}
\bibinfo{author}{Canuto, C.}, \bibinfo{author}{Hussaini, M.Y.},
  \bibinfo{author}{Quarteroni, A.}, \bibinfo{author}{Zang, T.A.},
  \bibinfo{year}{2006}.
\newblock \bibinfo{title}{Spectral methods}.
\newblock Scientific Computation, \bibinfo{publisher}{Springer-Verlag, Berlin}.
\bibitem[{Carstensen et~al.(2014)Carstensen, Feischl, Page and
  Praetorius}]{carstensen2014}
\bibinfo{author}{Carstensen, C.}, \bibinfo{author}{Feischl, M.},
  \bibinfo{author}{Page, M.}, \bibinfo{author}{Praetorius, D.},
  \bibinfo{year}{2014}.
\newblock \bibinfo{title}{Axioms of adaptivity}.
\newblock \bibinfo{journal}{Computers \& Mathematics with Applications}
  \bibinfo{volume}{67}, \bibinfo{pages}{1195--1253}.
\bibitem[{Carstensen and Merdon(2010)}]{carstensen2010}
\bibinfo{author}{Carstensen, C.}, \bibinfo{author}{Merdon, C.},
  \bibinfo{year}{2010}.
\newblock \bibinfo{title}{Estimator competition for {P}oisson problems}.
\newblock \bibinfo{journal}{Journal of Computational Mathematics}
  \bibinfo{volume}{28}, \bibinfo{pages}{309--330}.
\bibitem[{Cascon et~al.(2008)Cascon, Kreuzer, Nochetto and
  Siebert}]{cascon2008}
\bibinfo{author}{Cascon, J.M.}, \bibinfo{author}{Kreuzer, C.},
  \bibinfo{author}{Nochetto, R.H.}, \bibinfo{author}{Siebert, K.G.},
  \bibinfo{year}{2008}.
\newblock \bibinfo{title}{Quasi-optimal convergence rate for an adaptive finite
  element method}.
\newblock \bibinfo{journal}{SIAM Journal on Numerical Analysis}
  \bibinfo{volume}{46}, \bibinfo{pages}{2524--2550}.
\bibitem[{Chamoin and Legoll(2023)}]{chamoin2023}
\bibinfo{author}{Chamoin, L.}, \bibinfo{author}{Legoll, F.},
  \bibinfo{year}{2023}.
\newblock \bibinfo{title}{An introductory review on a posteriori error
  estimation in finite element computations}.
\newblock \bibinfo{journal}{SIAM Review} \bibinfo{volume}{65},
  \bibinfo{pages}{963--1028}.
\bibitem[{Chapelier et~al.(2024)Chapelier, Basile, Naddei, Plata, Couaillier
  and Laraufie}]{chapelier2024hp}
\bibinfo{author}{Chapelier, J.B.}, \bibinfo{author}{Basile, F.},
  \bibinfo{author}{Naddei, F.}, \bibinfo{author}{Plata, M.D.L.L.},
  \bibinfo{author}{Couaillier, V.}, \bibinfo{author}{Laraufie, R.},
  \bibinfo{year}{2024}.
\newblock \bibinfo{title}{hp adaptive discontinuous galerkin strategies driven
  by a posteriori error estimation with application to aeronautical flow
  problems}.
\newblock \bibinfo{journal}{Advances in Applied Mechanics (AAMS)}
  \bibinfo{volume}{58}.
\bibitem[{Chen et~al.(2024)Chen, Wang, Lian, Ma, Meng, Li, Ding and
  Bordas}]{chen2024reduced}
\bibinfo{author}{Chen, L.}, \bibinfo{author}{Wang, Z.}, \bibinfo{author}{Lian,
  H.}, \bibinfo{author}{Ma, Y.}, \bibinfo{author}{Meng, Z.},
  \bibinfo{author}{Li, P.}, \bibinfo{author}{Ding, C.},
  \bibinfo{author}{Bordas, S.P.}, \bibinfo{year}{2024}.
\newblock \bibinfo{title}{Reduced order isogeometric boundary element methods
  for cad-integrated shape optimization in electromagnetic scattering}.
\newblock \bibinfo{journal}{Computer Methods in Applied Mechanics and
  Engineering} \bibinfo{volume}{419}, \bibinfo{pages}{116654}.
\bibitem[{Ciarlet et~al.(2023)Ciarlet, Do and Madiot}]{ciarletmadioneutron}
\bibinfo{author}{Ciarlet, P.}, \bibinfo{author}{Do, M.H.},
  \bibinfo{author}{Madiot, F.}, \bibinfo{year}{2023}.
\newblock \bibinfo{title}{A \emph{posteriori} error estimates for mixed finite
  element discretizations of the neutron diffusion equations}.
\newblock \bibinfo{journal}{European Series in Applied and Industrial
  Mathematics (ESAIM): Mathematical Modelling and Numerical Analysis}
  \bibinfo{volume}{57}, \bibinfo{pages}{1--27}.
\bibitem[{Ciarlet(2002)}]{ciarlet2002}
\bibinfo{author}{Ciarlet, P.G.}, \bibinfo{year}{2002}.
\newblock \bibinfo{title}{The finite element method for elliptic problems}.
  volume~\bibinfo{volume}{40} of \textit{\bibinfo{series}{Classics in Applied
  Mathematics}}.
\newblock \bibinfo{publisher}{Society for Industrial and Applied Mathematics
  (SIAM), Philadelphia, PA}.
\bibitem[{Cicuttin et~al.(2021)Cicuttin, Ern and Pignet}]{cicuttin2021}
\bibinfo{author}{Cicuttin, M.}, \bibinfo{author}{Ern, A.},
  \bibinfo{author}{Pignet, N.}, \bibinfo{year}{2021}.
\newblock \bibinfo{title}{Hybrid high-order methods---a primer with
  applications to solid mechanics}.
\newblock SpringerBriefs in Mathematics, \bibinfo{publisher}{Springer, Cham}.
\bibitem[{Cockburn et~al.(2016)Cockburn, Di~Pietro and Ern}]{ern2016}
\bibinfo{author}{Cockburn, B.}, \bibinfo{author}{Di~Pietro, D.A.},
  \bibinfo{author}{Ern, A.}, \bibinfo{year}{2016}.
\newblock \bibinfo{title}{Bridging the hybrid high-order and hybridizable
  discontinuous {G}alerkin methods}.
\newblock \bibinfo{journal}{ESAIM. Mathematical Modelling and Numerical
  Analysis} \bibinfo{volume}{50}, \bibinfo{pages}{635--650}.
\bibitem[{Cohen et~al.(2001)Cohen, Dahmen and DeVore}]{cohen2001}
\bibinfo{author}{Cohen, A.}, \bibinfo{author}{Dahmen, W.},
  \bibinfo{author}{DeVore, R.}, \bibinfo{year}{2001}.
\newblock \bibinfo{title}{Adaptive wavelet methods for elliptic operator
  equations: convergence rates}.
\newblock \bibinfo{journal}{Mathematics of Computation} \bibinfo{volume}{70},
  \bibinfo{pages}{27--75}.
\bibitem[{Cohen and Masson(2000)}]{cohen2000}
\bibinfo{author}{Cohen, A.}, \bibinfo{author}{Masson, R.},
  \bibinfo{year}{2000}.
\newblock \bibinfo{title}{Wavelet adaptive method for second order elliptic
  problems: {Boundary} conditions and domain decomposition}.
\newblock \bibinfo{journal}{Numerische Mathematik} \bibinfo{volume}{86},
  \bibinfo{pages}{193--238}.
\bibitem[{{Conjungo Taumhas, Y. } et~al.(2023){Conjungo Taumhas, Y. },
  {Labeurthre, D.}, {Madiot, F.}, {Mula, O.} and {Taddei,
  T.}}]{taumhas2022impact}
\bibinfo{author}{{Conjungo Taumhas, Y. }}, \bibinfo{author}{{Labeurthre, D.}},
  \bibinfo{author}{{Madiot, F.}}, \bibinfo{author}{{Mula, O.}},
  \bibinfo{author}{{Taddei, T.}}, \bibinfo{year}{2023}.
\newblock \bibinfo{title}{Impact of physical model error on state estimation
  for neutronics applications*}.
\newblock \bibinfo{journal}{ESAIM: Proceedings} \bibinfo{volume}{73},
  \bibinfo{pages}{158--172}.
\bibitem[{Cottrell et~al.(2009)Cottrell, Hughes and Bazilevs}]{cottrell2009}
\bibinfo{author}{Cottrell, J.A.}, \bibinfo{author}{Hughes, T.J.R.},
  \bibinfo{author}{Bazilevs, Y.}, \bibinfo{year}{2009}.
\newblock \bibinfo{title}{Isogeometric analysis}.
\newblock \bibinfo{publisher}{John Wiley \& Sons, Ltd., Chichester}.
\bibitem[{Dakroub et~al.(2023)Dakroub, Faddoul, Omnes and Sayah}]{dakroub2023}
\bibinfo{author}{Dakroub, J.}, \bibinfo{author}{Faddoul, J.},
  \bibinfo{author}{Omnes, P.}, \bibinfo{author}{Sayah, T.},
  \bibinfo{year}{2023}.
\newblock \bibinfo{title}{A posteriori error estimates for the time-dependent
  {Navier}-{Stokes} system coupled with the convection-diffusion-reaction
  equation}.
\newblock \bibinfo{journal}{Advances in Computational Mathematics}
  \bibinfo{volume}{49}, \bibinfo{pages}{60}.
\newblock \bibinfo{note}{Id/No 67}.
\bibitem[{Davydov and Oanh(2011)}]{meshless_davydov2011adaptive}
\bibinfo{author}{Davydov, O.}, \bibinfo{author}{Oanh, D.T.},
  \bibinfo{year}{2011}.
\newblock \bibinfo{title}{Adaptive meshless centres and rbf stencils for
  poisson equation}.
\newblock \bibinfo{journal}{Journal of Computational Physics}
  \bibinfo{volume}{230}, \bibinfo{pages}{287--304}.
\bibitem[{Demkowicz and Gopalakrishnan(2010)}]{demkowicz2010}
\bibinfo{author}{Demkowicz, L.}, \bibinfo{author}{Gopalakrishnan, J.},
  \bibinfo{year}{2010}.
\newblock \bibinfo{title}{A class of discontinuous {P}etrov-{G}alerkin methods.
  {P}art {I}: the transport equation}.
\newblock \bibinfo{journal}{Computer Methods in Applied Mechanics and
  Engineering} \bibinfo{volume}{199}, \bibinfo{pages}{1558--1572}.
\bibitem[{Demkowicz et~al.(2012)Demkowicz, Gopalakrishnan and
  Niemi}]{demkowicz2012}
\bibinfo{author}{Demkowicz, L.}, \bibinfo{author}{Gopalakrishnan, J.},
  \bibinfo{author}{Niemi, A.H.}, \bibinfo{year}{2012}.
\newblock \bibinfo{title}{A class of discontinuous {P}etrov-{G}alerkin methods.
  {P}art {III}: {A}daptivity}.
\newblock \bibinfo{journal}{Applied Numerical Mathematics}
  \bibinfo{volume}{62}, \bibinfo{pages}{396--427}.
\bibitem[{Deshpande et~al.(2022)Deshpande, Lengiewicz and
  Bordas}]{deshpande2022probabilistic}
\bibinfo{author}{Deshpande, S.}, \bibinfo{author}{Lengiewicz, J.},
  \bibinfo{author}{Bordas, S.P.}, \bibinfo{year}{2022}.
\newblock \bibinfo{title}{Probabilistic deep learning for real-time large
  deformation simulations}.
\newblock \bibinfo{journal}{Computer Methods in Applied Mechanics and
  Engineering} \bibinfo{volume}{398}, \bibinfo{pages}{115307}.
\bibitem[{Deshpande et~al.(2023)Deshpande, Sosa, Bordas and
  Lengiewicz}]{deshpande2023convolution}
\bibinfo{author}{Deshpande, S.}, \bibinfo{author}{Sosa, R.I.},
  \bibinfo{author}{Bordas, S.}, \bibinfo{author}{Lengiewicz, J.},
  \bibinfo{year}{2023}.
\newblock \bibinfo{title}{Convolution, aggregation and attention based deep
  neural networks for accelerating simulations in mechanics}.
\newblock \bibinfo{journal}{Frontiers in Materials} \bibinfo{volume}{10},
  \bibinfo{pages}{1128954}.
\bibitem[{Di~Pietro and Droniou(2020)}]{dipietro2020}
\bibinfo{author}{Di~Pietro, D.A.}, \bibinfo{author}{Droniou, J.},
  \bibinfo{year}{2020}.
\newblock \bibinfo{title}{The hybrid high-order method for polytopal meshes}.
  volume~\bibinfo{volume}{19} of \textit{\bibinfo{series}{MS\&A. Modeling,
  Simulation and Applications}}.
\newblock \bibinfo{publisher}{Springer, Cham}.
\bibitem[{Di~Pietro and Ern(2012)}]{dipietro2012}
\bibinfo{author}{Di~Pietro, D.A.}, \bibinfo{author}{Ern, A.},
  \bibinfo{year}{2012}.
\newblock \bibinfo{title}{Mathematical aspects of discontinuous {G}alerkin
  methods}. volume~\bibinfo{volume}{69} of \textit{\bibinfo{series}{Mathematics
  \& Applications}}.
\newblock \bibinfo{publisher}{Springer, Heidelberg}.
\bibitem[{Dolbow and Belytschko(1998)}]{meshless_dolbow1998introduction}
\bibinfo{author}{Dolbow, J.}, \bibinfo{author}{Belytschko, T.},
  \bibinfo{year}{1998}.
\newblock \bibinfo{title}{{An introduction to programming the meshless Element
  FreeGalerkin method}}.
\newblock \bibinfo{journal}{Archives of Computational Methods in Engineering}
  \bibinfo{volume}{5}, \bibinfo{pages}{207--241}.
\bibitem[{Dong and Ern(2022)}]{dong2022}
\bibinfo{author}{Dong, Z.}, \bibinfo{author}{Ern, A.}, \bibinfo{year}{2022}.
\newblock \bibinfo{title}{Hybrid high-order and weak {Galerkin} methods for the
  biharmonic problem}.
\newblock \bibinfo{journal}{SIAM Journal on Numerical Analysis}
  \bibinfo{volume}{60}, \bibinfo{pages}{2626--2656}.
\bibitem[{D\"{o}rfler(1996)}]{dorfler1996}
\bibinfo{author}{D\"{o}rfler, W.}, \bibinfo{year}{1996}.
\newblock \bibinfo{title}{A convergent adaptive algorithm for {P}oisson's
  equation}.
\newblock \bibinfo{journal}{SIAM Journal on Numerical Analysis}
  \bibinfo{volume}{33}, \bibinfo{pages}{1106--1124}.
\bibitem[{Duflot and Bordas(2008)}]{duflot2008posteriori}
\bibinfo{author}{Duflot, M.}, \bibinfo{author}{Bordas, S.},
  \bibinfo{year}{2008}.
\newblock \bibinfo{title}{A posteriori error estimation for extended finite
  elements by an extended global recovery}.
\newblock \bibinfo{journal}{International Journal for Numerical Methods in
  Engineering} \bibinfo{volume}{76}, \bibinfo{pages}{1123--1138}.
\bibitem[{Duprez et~al.(2020)Duprez, Bordas, Bucki, Bui, Chouly, Lleras, Lobos,
  Lozinski, Rohan and Tomar}]{duprez2020soft}
\bibinfo{author}{Duprez, M.}, \bibinfo{author}{Bordas, S.P.A.},
  \bibinfo{author}{Bucki, M.}, \bibinfo{author}{Bui, H.P.},
  \bibinfo{author}{Chouly, F.}, \bibinfo{author}{Lleras, V.},
  \bibinfo{author}{Lobos, C.}, \bibinfo{author}{Lozinski, A.},
  \bibinfo{author}{Rohan, P.Y.}, \bibinfo{author}{Tomar, S.},
  \bibinfo{year}{2020}.
\newblock \bibinfo{title}{Quantifying discretization errors for soft tissue
  simulation in computer assisted surgery: a preliminary study}.
\newblock \bibinfo{journal}{Applied Mathematical Modelling. Simulation and
  Computation for Engineering and Environmental Systems} \bibinfo{volume}{77},
  \bibinfo{pages}{709--723}.
\bibitem[{Duprez and Lozinski(2020)}]{duprez2020fifem}
\bibinfo{author}{Duprez, M.}, \bibinfo{author}{Lozinski, A.},
  \bibinfo{year}{2020}.
\newblock \bibinfo{title}{{{\(\phi\)}}-{FEM}: a finite element method on
  domains defined by level-sets}.
\newblock \bibinfo{journal}{SIAM Journal on Numerical Analysis}
  \bibinfo{volume}{58}, \bibinfo{pages}{1008--1028}.
\bibitem[{Dur\'{a}n et~al.(1991)Dur\'{a}n, Muschietti and
  Rodr\'{\i}guez}]{duran1991}
\bibinfo{author}{Dur\'{a}n, R.}, \bibinfo{author}{Muschietti, M.A.},
  \bibinfo{author}{Rodr\'{\i}guez, R.}, \bibinfo{year}{1991}.
\newblock \bibinfo{title}{On the asymptotic exactness of error estimators for
  linear triangular finite elements}.
\newblock \bibinfo{journal}{Numerische Mathematik} \bibinfo{volume}{59},
  \bibinfo{pages}{107--127}.
\bibitem[{Eftimie et~al.(2023)Eftimie, Mavrodin and
  Bordas}]{eftimie2023digital}
\bibinfo{author}{Eftimie, R.}, \bibinfo{author}{Mavrodin, A.},
  \bibinfo{author}{Bordas, S.P.}, \bibinfo{year}{2023}.
\newblock \bibinfo{title}{From digital control to digital twins in medicine: A
  brief review and future perspectives}.
\newblock \bibinfo{journal}{Advances in Applied Mechanics}
  \bibinfo{volume}{56}, \bibinfo{pages}{323--368}.
\bibitem[{Eigel and Merdon(2016)}]{eigelI}
\bibinfo{author}{Eigel, M.}, \bibinfo{author}{Merdon, C.},
  \bibinfo{year}{2016}.
\newblock \bibinfo{title}{Local equilibration error estimators for guaranteed
  error control in adaptive stochastic higher-order {G}alerkin finite element
  methods}.
\newblock \bibinfo{journal}{SIAM/ASA Journal on Uncertainty Quantification}
  \bibinfo{volume}{4}, \bibinfo{pages}{1372--1397}.
\bibitem[{Eigel et~al.(2016)Eigel, Merdon and Neumann}]{eigelII}
\bibinfo{author}{Eigel, M.}, \bibinfo{author}{Merdon, C.},
  \bibinfo{author}{Neumann, J.}, \bibinfo{year}{2016}.
\newblock \bibinfo{title}{An adaptive multilevel {M}onte {C}arlo method with
  stochastic bounds for quantities of interest with uncertain data}.
\newblock \bibinfo{journal}{SIAM/ASA Journal on Uncertainty Quantification}
  \bibinfo{volume}{4}, \bibinfo{pages}{1219--1245}.
\bibitem[{Eriksson and Johnson(1985)}]{eriksson1985}
\bibinfo{author}{Eriksson, K.}, \bibinfo{author}{Johnson, C.},
  \bibinfo{year}{1985}.
\newblock \bibinfo{title}{Error estimates and automatic time step control for
  non-linear parabolic problems. {I}}.
\newblock \bibinfo{journal}{Bericht. Universit{\"a}t Jyv{\"a}skyl{\"a}.
  Mathematisches Institut} \bibinfo{volume}{31}.
\bibitem[{Ern and Guermond(2004)}]{ern-guermond-04}
\bibinfo{author}{Ern, A.}, \bibinfo{author}{Guermond, J.L.},
  \bibinfo{year}{2004}.
\newblock \bibinfo{title}{Theory and practice of finite elements}. volume
  \bibinfo{volume}{159} of \textit{\bibinfo{series}{Applied Mathematical
  Sciences}}.
\newblock \bibinfo{publisher}{Springer-Verlag}, \bibinfo{address}{New York}.
\bibitem[{Ern and Guermond(2021)}]{ern2021a}
\bibinfo{author}{Ern, A.}, \bibinfo{author}{Guermond, J.L.},
  \bibinfo{year}{2021}.
\newblock \bibinfo{title}{Finite elements. {I}---{A}pproximation and
  interpolation}. volume~\bibinfo{volume}{72} of \textit{\bibinfo{series}{Texts
  in Applied Mathematics}}.
\newblock \bibinfo{publisher}{Springer, Cham}.
\bibitem[{Ern et~al.(2017)Ern, Smears and Vohral\'{\i}k}]{Ern17}
\bibinfo{author}{Ern, A.}, \bibinfo{author}{Smears, I.},
  \bibinfo{author}{Vohral\'{\i}k, M.}, \bibinfo{year}{2017}.
\newblock \bibinfo{title}{Guaranteed, locally space-time efficient, and
  polynomial-degree robust a~posteriori error estimates for high-order
  discretizations of parabolic problems}.
\newblock \bibinfo{journal}{SIAM Journal on Numerical Analysis}
  \bibinfo{volume}{55}, \bibinfo{pages}{2811--2834}.
\bibitem[{Ern and Vohral\'{\i}k(2013)}]{ern2013}
\bibinfo{author}{Ern, A.}, \bibinfo{author}{Vohral\'{\i}k, M.},
  \bibinfo{year}{2013}.
\newblock \bibinfo{title}{Adaptive inexact {N}ewton methods with a posteriori
  stopping criteria for nonlinear diffusion {PDE}s}.
\newblock \bibinfo{journal}{SIAM Journal on Scientific Computing}
  \bibinfo{volume}{35}, \bibinfo{pages}{A1761--A1791}.
\bibitem[{Feischl(2019)}]{Feischl19}
\bibinfo{author}{Feischl, M.}, \bibinfo{year}{2019}.
\newblock \bibinfo{title}{Optimality of a standard adaptive finite element
  method for the {S}tokes problem}.
\newblock \bibinfo{journal}{SIAM Journal on Numerical Analysis}
  \bibinfo{volume}{57}, \bibinfo{pages}{1124--1157}.
\bibitem[{Feischl(2022)}]{Feischl22}
\bibinfo{author}{Feischl, M.}, \bibinfo{year}{2022}.
\newblock \bibinfo{title}{Inf-sup stability implies quasi-orthogonality}.
\newblock \bibinfo{journal}{Mathematics of Computation} \bibinfo{volume}{91},
  \bibinfo{pages}{2059--2094}.
\bibitem[{Feischl et~al.(2014)Feischl, F\"{u}hrer and
  Praetorius}]{FeischlFuhrerPraetorius14}
\bibinfo{author}{Feischl, M.}, \bibinfo{author}{F\"{u}hrer, T.},
  \bibinfo{author}{Praetorius, D.}, \bibinfo{year}{2014}.
\newblock \bibinfo{title}{Adaptive {FEM} with optimal convergence rates for a
  certain class of nonsymmetric and possibly nonlinear problems}.
\newblock \bibinfo{journal}{SIAM Journal on Numerical Analysis}
  \bibinfo{volume}{52}, \bibinfo{pages}{601--625}.
\bibitem[{Feischl et~al.(2016)Feischl, Praetorius and van~der
  Zee}]{FeischlPraetoriusZee16}
\bibinfo{author}{Feischl, M.}, \bibinfo{author}{Praetorius, D.},
  \bibinfo{author}{van~der Zee, K.G.}, \bibinfo{year}{2016}.
\newblock \bibinfo{title}{An abstract analysis of optimal goal-oriented
  adaptivity}.
\newblock \bibinfo{journal}{SIAM Journal on Numerical Analysis}
  \bibinfo{volume}{54}, \bibinfo{pages}{1423--1448}.
\bibitem[{Ferraz-Leite et~al.(2010)Ferraz-Leite, Ortner and
  Praetorius}]{Ferraz-LeiteOrtnerPraetorius10}
\bibinfo{author}{Ferraz-Leite, S.}, \bibinfo{author}{Ortner, C.},
  \bibinfo{author}{Praetorius, D.}, \bibinfo{year}{2010}.
\newblock \bibinfo{title}{Convergence of simple adaptive {G}alerkin schemes
  based on {$h-h/2$} error estimators}.
\newblock \bibinfo{journal}{Numerische Mathematik} \bibinfo{volume}{116},
  \bibinfo{pages}{291--316}.
\bibitem[{F\'{e}votte et~al.(2024)F\'{e}votte, Rappaport and
  Vohral\'{\i}k}]{Vohralik24}
\bibinfo{author}{F\'{e}votte, F.}, \bibinfo{author}{Rappaport, A.},
  \bibinfo{author}{Vohral\'{\i}k, M.}, \bibinfo{year}{2024}.
\newblock \bibinfo{title}{Adaptive regularization, discretization, and
  linearization for nonsmooth problems based on primal-dual gap estimators}.
\newblock \bibinfo{journal}{Computer Methods in Applied Mechanics and
  Engineering} \bibinfo{volume}{418}, \bibinfo{pages}{Paper No. 116558, 33}.
\bibitem[{Fortin(2024)}]{fortin2023anisotropic}
\bibinfo{author}{Fortin, A.}, \bibinfo{year}{2024}.
\newblock \bibinfo{title}{An anisotropic mesh adaptation method based on
  gradient recovery and optimal shape elements}.
\newblock \bibinfo{journal}{Advances in Applied Mechanics (AAMS)}
  \bibinfo{volume}{58}.
\bibitem[{Gander and Wanner(2012)}]{gander2012}
\bibinfo{author}{Gander, M.J.}, \bibinfo{author}{Wanner, G.},
  \bibinfo{year}{2012}.
\newblock \bibinfo{title}{From {Euler}, {Ritz}, and {Galerkin} to modern
  computing}.
\newblock \bibinfo{journal}{SIAM Review} \bibinfo{volume}{54},
  \bibinfo{pages}{627--666}.
\bibitem[{Gerasimov et~al.(2012)Gerasimov, R{\"u}ter and Stein}]{gerasimov2012}
\bibinfo{author}{Gerasimov, T.}, \bibinfo{author}{R{\"u}ter, M.},
  \bibinfo{author}{Stein, E.}, \bibinfo{year}{2012}.
\newblock \bibinfo{title}{An explicit residual-type error estimator for
  {Q1}-quadrilateral extended finite element method in two-dimensional linear
  elastic fracture mechanics}.
\newblock \bibinfo{journal}{International Journal for Numerical Methods in
  Engineering} \bibinfo{volume}{90}, \bibinfo{pages}{1118--1155}.
\bibitem[{Giles and S{\"u}li(2002)}]{giles2002}
\bibinfo{author}{Giles, M.B.}, \bibinfo{author}{S{\"u}li, E.},
  \bibinfo{year}{2002}.
\newblock \bibinfo{title}{Adjoint methods for {PDEs}: a posteriori error
  analysis and postprocessing by duality}.
\newblock \bibinfo{journal}{Acta Numerica} \bibinfo{volume}{11},
  \bibinfo{pages}{145--236}.
\bibitem[{Glowinski et~al.(1994)Glowinski, Pan and
  P{\'e}riaux}]{glowinski-pan-periaux-94}
\bibinfo{author}{Glowinski, R.}, \bibinfo{author}{Pan, T.},
  \bibinfo{author}{P{\'e}riaux, J.}, \bibinfo{year}{1994}.
\newblock \bibinfo{title}{A fictitious domain method for {D}irichlet problem
  and applications}.
\newblock \bibinfo{journal}{Computer Methods in Applied Mechanics and
  Engineering} \bibinfo{volume}{111}, \bibinfo{pages}{283--303}.
\bibitem[{Gonz{\'a}lez-Estrada et~al.(2014)Gonz{\'a}lez-Estrada, Nadal,
  R{\'o}denas, Kerfriden, Bordas and Fuenmayor}]{gonzalez2014}
\bibinfo{author}{Gonz{\'a}lez-Estrada, O.A.}, \bibinfo{author}{Nadal, E.},
  \bibinfo{author}{R{\'o}denas, J.J.}, \bibinfo{author}{Kerfriden, P.},
  \bibinfo{author}{Bordas, S.P.A.}, \bibinfo{author}{Fuenmayor, F.J.},
  \bibinfo{year}{2014}.
\newblock \bibinfo{title}{Mesh adaptivity driven by goal-oriented locally
  equilibrated superconvergent patch recovery}.
\newblock \bibinfo{journal}{Computational Mechanics} \bibinfo{volume}{53},
  \bibinfo{pages}{957--976}.
\bibitem[{Gonz{\'a}lez-Estrada et~al.(2021)Gonz{\'a}lez-Estrada, Natarajan,
  R{\'o}denas and Bordas}]{gonzalez2021error}
\bibinfo{author}{Gonz{\'a}lez-Estrada, O.A.}, \bibinfo{author}{Natarajan, S.},
  \bibinfo{author}{R{\'o}denas, J.J.}, \bibinfo{author}{Bordas, S.P.},
  \bibinfo{year}{2021}.
\newblock \bibinfo{title}{Error estimation for the polygonal finite element
  method for smooth and singular linear elasticity}.
\newblock \bibinfo{journal}{Computers \& Mathematics with Applications}
  \bibinfo{volume}{92}, \bibinfo{pages}{109--119}.
\bibitem[{Gonz{\'a}lez-Estrada et~al.(2013)Gonz{\'a}lez-Estrada, Natarajan,
  R{\'o}denas, Nguyen-Xuan and Bordas}]{gonzalez2013efficient}
\bibinfo{author}{Gonz{\'a}lez-Estrada, O.A.}, \bibinfo{author}{Natarajan, S.},
  \bibinfo{author}{R{\'o}denas, J.J.}, \bibinfo{author}{Nguyen-Xuan, H.},
  \bibinfo{author}{Bordas, S.P.}, \bibinfo{year}{2013}.
\newblock \bibinfo{title}{Efficient recovery-based error estimation for the
  smoothed finite element method for smooth and singular linear elasticity}.
\newblock \bibinfo{journal}{Computational Mechanics} \bibinfo{volume}{52},
  \bibinfo{pages}{37--52}.
\bibitem[{Gonz{\'a}lez-Estrada et~al.(2012)Gonz{\'a}lez-Estrada, R{\'o}denas,
  Bordas, Duflot, Kerfriden and Giner}]{gonzalez2012role}
\bibinfo{author}{Gonz{\'a}lez-Estrada, O.A.}, \bibinfo{author}{R{\'o}denas,
  J.J.}, \bibinfo{author}{Bordas, S.P.A.}, \bibinfo{author}{Duflot, M.},
  \bibinfo{author}{Kerfriden, P.}, \bibinfo{author}{Giner, E.},
  \bibinfo{year}{2012}.
\newblock \bibinfo{title}{On the role of enrichment and statical admissibility
  of recovered fields in a posteriori error estimation for enriched finite
  element methods}.
\newblock \bibinfo{journal}{Engineering Computations} \bibinfo{volume}{29},
  \bibinfo{pages}{814--841}.
\bibitem[{Gonz{\'a}lez-Estrada et~al.(2023)Gonz{\'a}lez-Estrada,
  R{\'o}denas~Garc{\'\i}a and Bordas}]{gonzalez2023recovery}
\bibinfo{author}{Gonz{\'a}lez-Estrada, O.A.},
  \bibinfo{author}{R{\'o}denas~Garc{\'\i}a, J.J.}, \bibinfo{author}{Bordas,
  S.P.}, \bibinfo{year}{2023}.
\newblock \bibinfo{title}{Recovery-based error estimation and bounding in
  xfem}.
\newblock \bibinfo{journal}{Partition of Unity Methods} ,
  \bibinfo{pages}{145--190}.
\bibitem[{Goury et~al.(2016)Goury, Amsallem, Bordas, Liu and
  Kerfriden}]{goury2016automatised}
\bibinfo{author}{Goury, O.}, \bibinfo{author}{Amsallem, D.},
  \bibinfo{author}{Bordas, S.P.A.}, \bibinfo{author}{Liu, W.K.},
  \bibinfo{author}{Kerfriden, P.}, \bibinfo{year}{2016}.
\newblock \bibinfo{title}{Automatised selection of load paths to construct
  reduced-order models in computational damage micromechanics: from
  dissipation-driven random selection to bayesian optimization}.
\newblock \bibinfo{journal}{Computational Mechanics} \bibinfo{volume}{58},
  \bibinfo{pages}{213--234}.
\bibitem[{Grepl et~al.(2007)Grepl, Maday, Nguyen and Patera}]{grepl2007}
\bibinfo{author}{Grepl, M.A.}, \bibinfo{author}{Maday, Y.},
  \bibinfo{author}{Nguyen, N.C.}, \bibinfo{author}{Patera, A.T.},
  \bibinfo{year}{2007}.
\newblock \bibinfo{title}{Efficient reduced-basis treatment of nonaffine and
  nonlinear partial differential equations}.
\newblock \bibinfo{journal}{M2AN. Mathematical Modelling and Numerical
  Analysis} \bibinfo{volume}{41}, \bibinfo{pages}{575--605}.
\bibitem[{Guignard and Nobile(2018)}]{guignardII}
\bibinfo{author}{Guignard, D.}, \bibinfo{author}{Nobile, F.},
  \bibinfo{year}{2018}.
\newblock \bibinfo{title}{A posteriori error estimation for the stochastic
  collocation finite element method}.
\newblock \bibinfo{journal}{SIAM Journal on Numerical Analysis}
  \bibinfo{volume}{56}, \bibinfo{pages}{3121--3143}.
\bibitem[{Guignard et~al.(2016)Guignard, Nobile and Picasso}]{guignardI}
\bibinfo{author}{Guignard, D.}, \bibinfo{author}{Nobile, F.},
  \bibinfo{author}{Picasso, M.}, \bibinfo{year}{2016}.
\newblock \bibinfo{title}{A posteriori error estimation for elliptic partial
  differential equations with small uncertainties}.
\newblock \bibinfo{journal}{Numerical Methods for Partial Differential
  Equations. An International Journal} \bibinfo{volume}{32},
  \bibinfo{pages}{175--212}.
\bibitem[{Gustafsson(2024)}]{gustafsson2024adaptive}
\bibinfo{author}{Gustafsson, T.}, \bibinfo{year}{2024}.
\newblock \bibinfo{title}{Adaptive finite elements for obstacle problems}.
\newblock \bibinfo{journal}{Advances in Applied Mechanics (AAMS)}
  \bibinfo{volume}{58}.
\bibitem[{Haberl et~al.(2021)Haberl, Praetorius, Schimanko and
  Vohral{\'{\i}}k}]{haberl2021}
\bibinfo{author}{Haberl, A.}, \bibinfo{author}{Praetorius, D.},
  \bibinfo{author}{Schimanko, S.}, \bibinfo{author}{Vohral{\'{\i}}k, M.},
  \bibinfo{year}{2021}.
\newblock \bibinfo{title}{Convergence and quasi-optimal cost of adaptive
  algorithms for nonlinear operators including iterative linearization and
  algebraic solver}.
\newblock \bibinfo{journal}{Numerische Mathematik} \bibinfo{volume}{147},
  \bibinfo{pages}{679--725}.
\bibitem[{Han(2005)}]{han2005book}
\bibinfo{author}{Han, W.}, \bibinfo{year}{2005}.
\newblock \bibinfo{title}{A posteriori error analysis via duality theory.
  {With} applications in modeling and numerical approximations}.
  volume~\bibinfo{volume}{8} of \textit{\bibinfo{series}{Advances in Mechanics
  and Mathematics}}.
\newblock \bibinfo{publisher}{New York, NY: Springer}.
\bibitem[{Haslinger and Renard(2009)}]{haslinger2009}
\bibinfo{author}{Haslinger, J.}, \bibinfo{author}{Renard, Y.},
  \bibinfo{year}{2009}.
\newblock \bibinfo{title}{A new fictitious domain approach inspired by the
  extended finite element method}.
\newblock \bibinfo{journal}{SIAM Journal on Numerical Analysis}
  \bibinfo{volume}{47}, \bibinfo{pages}{1474--1499}.
\bibitem[{Hauseux et~al.(2017a)Hauseux, Hale and
  Bordas}]{hauseux2017accelerating}
\bibinfo{author}{Hauseux, P.}, \bibinfo{author}{Hale, J.S.},
  \bibinfo{author}{Bordas, S.P.}, \bibinfo{year}{2017}a.
\newblock \bibinfo{title}{Accelerating monte carlo estimation with derivatives
  of high-level finite element models}.
\newblock \bibinfo{journal}{Computer Methods in Applied Mechanics and
  Engineering} \bibinfo{volume}{318}, \bibinfo{pages}{917--936}.
\bibitem[{Hauseux et~al.(2017b)Hauseux, Hale and
  Bordas}]{hauseux2017calculating}
\bibinfo{author}{Hauseux, P.}, \bibinfo{author}{Hale, J.S.},
  \bibinfo{author}{Bordas, S.P.}, \bibinfo{year}{2017}b.
\newblock \bibinfo{title}{Calculating the malliavin derivative of some
  stochastic mechanics problems}.
\newblock \bibinfo{journal}{Plos one} \bibinfo{volume}{12},
  \bibinfo{pages}{e0189994}.
\bibitem[{Hauseux et~al.(2018)Hauseux, Hale, Cotin and
  Bordas}]{hauseux2018quantifying}
\bibinfo{author}{Hauseux, P.}, \bibinfo{author}{Hale, J.S.},
  \bibinfo{author}{Cotin, S.}, \bibinfo{author}{Bordas, S.P.},
  \bibinfo{year}{2018}.
\newblock \bibinfo{title}{Quantifying the uncertainty in a hyperelastic soft
  tissue model with stochastic parameters}.
\newblock \bibinfo{journal}{Applied Mathematical Modelling}
  \bibinfo{volume}{62}, \bibinfo{pages}{86--102}.
\bibitem[{Hesthaven et~al.(2016)Hesthaven, Rozza and Stamm}]{hesthaven2016}
\bibinfo{author}{Hesthaven, J.S.}, \bibinfo{author}{Rozza, G.},
  \bibinfo{author}{Stamm, B.}, \bibinfo{year}{2016}.
\newblock \bibinfo{title}{Certified reduced basis methods for parametrized
  partial differential equations}.
\newblock SpringerBriefs in Mathematics, \bibinfo{publisher}{Springer, Cham;
  BCAM Basque Center for Applied Mathematics, Bilbao}.
\newblock \bibinfo{note}{BCAM SpringerBriefs}.
\bibitem[{Hild et~al.(2009)Hild, Lleras and Renard}]{hild2009}
\bibinfo{author}{Hild, P.}, \bibinfo{author}{Lleras, V.},
  \bibinfo{author}{Renard, Y.}, \bibinfo{year}{2009}.
\newblock \bibinfo{title}{A posteriori error analysis for {Poisson}'s equation
  approximated by {XFEM}}.
\newblock \bibinfo{journal}{European Series in Applied and Industrial
  Mathematics (ESAIM): Proceedings} \bibinfo{volume}{27},
  \bibinfo{pages}{107--121}.
\bibitem[{Hoang et~al.(2021)Hoang, Choi and Carlberg}]{hoang2021domain}
\bibinfo{author}{Hoang, C.}, \bibinfo{author}{Choi, Y.},
  \bibinfo{author}{Carlberg, K.}, \bibinfo{year}{2021}.
\newblock \bibinfo{title}{{Domain-decomposition least-squares Petrov--Galerkin
  (DD-LSPG) nonlinear model reduction}}.
\newblock \bibinfo{journal}{Computer Methods in Applied Mechanics and
  Engineering} \bibinfo{volume}{384}, \bibinfo{pages}{113997}.
\bibitem[{Hoang et~al.(2022)Hoang, Chowdhary, Lee and
  Ray}]{hoang2022projection}
\bibinfo{author}{Hoang, C.}, \bibinfo{author}{Chowdhary, K.},
  \bibinfo{author}{Lee, K.}, \bibinfo{author}{Ray, J.}, \bibinfo{year}{2022}.
\newblock \bibinfo{title}{Projection-based model reduction of dynamical systems
  using space--time subspace and machine learning}.
\newblock \bibinfo{journal}{Computer Methods in Applied Mechanics and
  Engineering} \bibinfo{volume}{389}, \bibinfo{pages}{114341}.
\bibitem[{Hoang et~al.(2016)Hoang, Kerfriden and Bordas}]{hoang2016fast}
\bibinfo{author}{Hoang, K.}, \bibinfo{author}{Kerfriden, P.},
  \bibinfo{author}{Bordas, S.P.A.}, \bibinfo{year}{2016}.
\newblock \bibinfo{title}{A fast, certified and “tuning free” two-field
  reduced basis method for the metamodelling of affinely-parametrised
  elasticity problems}.
\newblock \bibinfo{journal}{Computer Methods in Applied Mechanics and
  Engineering} \bibinfo{volume}{298}, \bibinfo{pages}{121--158}.
\bibitem[{Hu et~al.(2018)Hu, Chouly, Hu, Cheng and Bordas}]{hu2018}
\bibinfo{author}{Hu, Q.}, \bibinfo{author}{Chouly, F.}, \bibinfo{author}{Hu,
  P.}, \bibinfo{author}{Cheng, G.}, \bibinfo{author}{Bordas, S.P.A.},
  \bibinfo{year}{2018}.
\newblock \bibinfo{title}{Skew-symmetric {Nitsche}'s formulation in
  isogeometric analysis: {Dirichlet} and symmetry conditions, patch coupling
  and frictionless contact}.
\newblock \bibinfo{journal}{Computer Methods in Applied Mechanics and
  Engineering} \bibinfo{volume}{341}, \bibinfo{pages}{188--220}.
\bibitem[{Jacquemin et~al.(2023)Jacquemin, Suchde and
  Bordas}]{jacquemin2023smart}
\bibinfo{author}{Jacquemin, T.}, \bibinfo{author}{Suchde, P.},
  \bibinfo{author}{Bordas, S.P.}, \bibinfo{year}{2023}.
\newblock \bibinfo{title}{Smart cloud collocation: geometry-aware adaptivity
  directly from cad}.
\newblock \bibinfo{journal}{Computer-Aided Design} \bibinfo{volume}{154},
  \bibinfo{pages}{103409}.
\bibitem[{Jansari et~al.(2022a)Jansari, Videla, Natarajan, Bordas and
  Atroshchenko}]{IGA_GIFT_jansari2022adaptive}
\bibinfo{author}{Jansari, C.}, \bibinfo{author}{Videla, J.},
  \bibinfo{author}{Natarajan, S.}, \bibinfo{author}{Bordas, S.P.},
  \bibinfo{author}{Atroshchenko, E.}, \bibinfo{year}{2022}a.
\newblock \bibinfo{title}{Adaptive enriched geometry independent field
  approximation for 2d time-harmonic acoustics}.
\newblock \bibinfo{journal}{Computers \& Structures} \bibinfo{volume}{263},
  \bibinfo{pages}{106728}.
\bibitem[{Jansari et~al.(2022b)Jansari, Videla, Natarajan, Bordas and
  Atroshchenko}]{jansari2022adaptive}
\bibinfo{author}{Jansari, C.}, \bibinfo{author}{Videla, J.},
  \bibinfo{author}{Natarajan, S.}, \bibinfo{author}{Bordas, S.P.},
  \bibinfo{author}{Atroshchenko, E.}, \bibinfo{year}{2022}b.
\newblock \bibinfo{title}{Adaptive enriched geometry independent field
  approximation for 2d time-harmonic acoustics}.
\newblock \bibinfo{journal}{Computers \& Structures} \bibinfo{volume}{263},
  \bibinfo{pages}{106728}.
\bibitem[{Jia et~al.(2019)Jia, Anitescu, Zhang and Rabczuk}]{jia2019adaptive}
\bibinfo{author}{Jia, Y.}, \bibinfo{author}{Anitescu, C.},
  \bibinfo{author}{Zhang, Y.J.}, \bibinfo{author}{Rabczuk, T.},
  \bibinfo{year}{2019}.
\newblock \bibinfo{title}{An adaptive isogeometric analysis collocation method
  with a recovery-based error estimator}.
\newblock \bibinfo{journal}{Computer Methods in Applied Mechanics and
  Engineering} \bibinfo{volume}{345}, \bibinfo{pages}{52--74}.
\bibitem[{Jin et~al.(2017)Jin, Gonz{\'a}lez-Estrada, Pierard and
  Bordas}]{jin2017error}
\bibinfo{author}{Jin, Y.}, \bibinfo{author}{Gonz{\'a}lez-Estrada, O.},
  \bibinfo{author}{Pierard, O.}, \bibinfo{author}{Bordas, S.},
  \bibinfo{year}{2017}.
\newblock \bibinfo{title}{Error-controlled adaptive extended finite element
  method for 3d linear elastic crack propagation}.
\newblock \bibinfo{journal}{Computer Methods in Applied Mechanics and
  Engineering} \bibinfo{volume}{318}, \bibinfo{pages}{319--348}.
\bibitem[{Karaivanov and Petrushev(2003)}]{karaivanov2003}
\bibinfo{author}{Karaivanov, B.}, \bibinfo{author}{Petrushev, P.},
  \bibinfo{year}{2003}.
\newblock \bibinfo{title}{Nonlinear piecewise polynomial approximation beyond
  {Besov} spaces}.
\newblock \bibinfo{journal}{Applied and Computational Harmonic Analysis}
  \bibinfo{volume}{15}, \bibinfo{pages}{177--223}.
\bibitem[{Kerfriden et~al.(2011)Kerfriden, Gosselet, Adhikari and
  Bordas}]{kerfriden2011}
\bibinfo{author}{Kerfriden, P.}, \bibinfo{author}{Gosselet, P.},
  \bibinfo{author}{Adhikari, S.}, \bibinfo{author}{Bordas, S.P.A.},
  \bibinfo{year}{2011}.
\newblock \bibinfo{title}{Bridging proper orthogonal decomposition methods and
  augmented {N}ewton-{K}rylov algorithms: an adaptive model order reduction for
  highly nonlinear mechanical problems}.
\newblock \bibinfo{journal}{Computer Methods in Applied Mechanics and
  Engineering} \bibinfo{volume}{200}, \bibinfo{pages}{850--866}.
\bibitem[{Kerfriden et~al.(2014)Kerfriden, R{\'o}denas and
  Bordas}]{kerfriden2014certification}
\bibinfo{author}{Kerfriden, P.}, \bibinfo{author}{R{\'o}denas, J.J.},
  \bibinfo{author}{Bordas, S.A.}, \bibinfo{year}{2014}.
\newblock \bibinfo{title}{Certification of projection-based reduced order
  modelling in computational homogenisation by the constitutive relation
  error}.
\newblock \bibinfo{journal}{International Journal for Numerical Methods in
  Engineering} \bibinfo{volume}{97}, \bibinfo{pages}{395--422}.
\bibitem[{Kreuzer and Georgoulis(2018)}]{Kreuzer18}
\bibinfo{author}{Kreuzer, C.}, \bibinfo{author}{Georgoulis, E.H.},
  \bibinfo{year}{2018}.
\newblock \bibinfo{title}{Convergence of adaptive discontinuous {G}alerkin
  methods}.
\newblock \bibinfo{journal}{Mathematics of Computation} \bibinfo{volume}{87},
  \bibinfo{pages}{2611--2640}.
\bibitem[{Ladeveze and Leguillon(1983)}]{ladeveze1983}
\bibinfo{author}{Ladeveze, P.}, \bibinfo{author}{Leguillon, D.},
  \bibinfo{year}{1983}.
\newblock \bibinfo{title}{Error estimate procedure in the finite element method
  and applications}.
\newblock \bibinfo{journal}{SIAM Journal on Numerical Analysis}
  \bibinfo{volume}{20}, \bibinfo{pages}{485--509}.
\bibitem[{Langer et~al.(2016)Langer, Matculevich and
  Repin}]{langer2016posteriori}
\bibinfo{author}{Langer, U.}, \bibinfo{author}{Matculevich, S.},
  \bibinfo{author}{Repin, S.}, \bibinfo{year}{2016}.
\newblock \bibinfo{title}{A posteriori error estimates for space-time iga
  approximations to parabolic initial boundary value problems}.
\newblock \bibinfo{journal}{arXiv preprint arXiv:1612.08998} .
\bibitem[{Langer et~al.(2019)Langer, Matculevich and
  Repin}]{langer2019guaranteed}
\bibinfo{author}{Langer, U.}, \bibinfo{author}{Matculevich, S.},
  \bibinfo{author}{Repin, S.}, \bibinfo{year}{2019}.
\newblock \bibinfo{title}{Guaranteed error bounds and local indicators for
  adaptive solvers using stabilised space--time iga approximations to parabolic
  problems}.
\newblock \bibinfo{journal}{Computers \& Mathematics with Applications}
  \bibinfo{volume}{78}, \bibinfo{pages}{2641--2671}.
\bibitem[{Lemaire(2021)}]{lemaire2021}
\bibinfo{author}{Lemaire, S.}, \bibinfo{year}{2021}.
\newblock \bibinfo{title}{Bridging the hybrid high-order and virtual element
  methods}.
\newblock \bibinfo{journal}{IMA Journal of Numerical Analysis}
  \bibinfo{volume}{41}, \bibinfo{pages}{549--593}.
\bibitem[{Lesaint and Raviart(1974)}]{lesaint1974}
\bibinfo{author}{Lesaint, P.}, \bibinfo{author}{Raviart, P.A.},
  \bibinfo{year}{1974}.
\newblock \bibinfo{title}{On a finite element method for solving the neutron
  transport equation}, in: \bibinfo{booktitle}{Mathematical aspects of finite
  elements in partial differential equations ({P}roc. {S}ympos., {M}ath. {R}es.
  {C}enter, {U}niv. {W}isconsin, {M}adison, {W}is., 1974)}.
  \bibinfo{publisher}{Academic Press, New York-London}, pp.
  \bibinfo{pages}{89--123}.
\bibitem[{Li et~al.(2020)Li, Nguyen-Thanh, Huang and
  Zhou}]{meshless_shells_adaptive_li2020adaptive}
\bibinfo{author}{Li, W.}, \bibinfo{author}{Nguyen-Thanh, N.},
  \bibinfo{author}{Huang, J.}, \bibinfo{author}{Zhou, K.},
  \bibinfo{year}{2020}.
\newblock \bibinfo{title}{Adaptive analysis of crack propagation in thin-shell
  structures via an isogeometric-meshfree moving least-squares approach}.
\newblock \bibinfo{journal}{Computer Methods in Applied Mechanics and
  Engineering} \bibinfo{volume}{358}, \bibinfo{pages}{112613}.
\bibitem[{Liu et~al.(2007)Liu, Dai and Nguyen}]{liu2007}
\bibinfo{author}{Liu, G.R.}, \bibinfo{author}{Dai, K.Y.},
  \bibinfo{author}{Nguyen, T.T.}, \bibinfo{year}{2007}.
\newblock \bibinfo{title}{A smoothed finite element method for mechanics
  problems}.
\newblock \bibinfo{journal}{Computational Mechanics} \bibinfo{volume}{39},
  \bibinfo{pages}{859--877}.
\bibitem[{Liu et~al.(2017)Liu, Rami{\`e}re and Lebon}]{liu2017}
\bibinfo{author}{Liu, H.}, \bibinfo{author}{Rami{\`e}re, I.},
  \bibinfo{author}{Lebon, F.}, \bibinfo{year}{2017}.
\newblock \bibinfo{title}{On the coupling of local multilevel mesh refinement
  and {ZZ} methods for unilateral frictional contact problems in
  elastostatics}.
\newblock \bibinfo{journal}{Computer Methods in Applied Mechanics and
  Engineering} \bibinfo{volume}{323}, \bibinfo{pages}{1--26}.
\bibitem[{Lozinski(2019)}]{lozinski2019poly}
\bibinfo{author}{Lozinski, A.}, \bibinfo{year}{2019}.
\newblock \bibinfo{title}{A primal discontinuous {G}alerkin method with static
  condensation on very general meshes}.
\newblock \bibinfo{journal}{Numerische Mathematik} \bibinfo{volume}{143},
  \bibinfo{pages}{583--604}.
\bibitem[{Matculevich(2018)}]{IGA_functional_matculevich2018functional}
\bibinfo{author}{Matculevich, S.}, \bibinfo{year}{2018}.
\newblock \bibinfo{title}{Functional approach to the error control in adaptive
  iga schemes for elliptic boundary value problems}.
\newblock \bibinfo{journal}{Journal of Computational and Applied Mathematics}
  \bibinfo{volume}{344}, \bibinfo{pages}{394--423}.
\bibitem[{Mo\"{e}s et~al.(2006)Mo\"{e}s, B\'{e}chet and Tourbier}]{moes2006}
\bibinfo{author}{Mo\"{e}s, N.}, \bibinfo{author}{B\'{e}chet, E.},
  \bibinfo{author}{Tourbier, M.}, \bibinfo{year}{2006}.
\newblock \bibinfo{title}{Imposing {D}irichlet boundary conditions in the
  extended finite element method}.
\newblock \bibinfo{journal}{International Journal for Numerical Methods in
  Engineering} \bibinfo{volume}{67}, \bibinfo{pages}{1641--1669}.
\bibitem[{Mommer and Stevenson(2009)}]{MommerStevenson09}
\bibinfo{author}{Mommer, M.S.}, \bibinfo{author}{Stevenson, R.},
  \bibinfo{year}{2009}.
\newblock \bibinfo{title}{A goal-oriented adaptive finite element method with
  convergence rates}.
\newblock \bibinfo{journal}{SIAM Journal on Numerical Analysis}
  \bibinfo{volume}{47}, \bibinfo{pages}{861--886}.
\bibitem[{Monasse and Perrier(1998)}]{monasse1998}
\bibinfo{author}{Monasse, P.}, \bibinfo{author}{Perrier, V.},
  \bibinfo{year}{1998}.
\newblock \bibinfo{title}{Orthonormal wavelet bases adapted for partial
  differential equations with boundary conditions}.
\newblock \bibinfo{journal}{SIAM Journal on Mathematical Analysis}
  \bibinfo{volume}{29}, \bibinfo{pages}{1040--1065}.
\bibitem[{Morin et~al.(2000)Morin, Nochetto and Siebert}]{morin2000}
\bibinfo{author}{Morin, P.}, \bibinfo{author}{Nochetto, R.H.},
  \bibinfo{author}{Siebert, K.G.}, \bibinfo{year}{2000}.
\newblock \bibinfo{title}{Data oscillation and convergence of adaptive {FEM}}.
\newblock \bibinfo{journal}{SIAM Journal on Numerical Analysis}
  \bibinfo{volume}{38}, \bibinfo{pages}{466--488}.
\bibitem[{Nassreddine et~al.(2022)Nassreddine, Omnes and
  Sayah}]{nassreddine2022}
\bibinfo{author}{Nassreddine, G.}, \bibinfo{author}{Omnes, P.},
  \bibinfo{author}{Sayah, T.}, \bibinfo{year}{2022}.
\newblock \bibinfo{title}{A posteriori error estimates for the large eddy
  simulation applied to stationary {Navier}-{Stokes} equations}.
\newblock \bibinfo{journal}{Numerical Methods for Partial Differential
  Equations} \bibinfo{volume}{38}, \bibinfo{pages}{1468--1498}.
\bibitem[{Nassreddine et~al.(2023)Nassreddine, Omnes and
  Sayah}]{nassreddine2023}
\bibinfo{author}{Nassreddine, G.}, \bibinfo{author}{Omnes, P.},
  \bibinfo{author}{Sayah, T.}, \bibinfo{year}{2023}.
\newblock \bibinfo{title}{\emph{A posteriori} error estimates for the large
  eddy simulation applied to incompressible fluids}.
\newblock \bibinfo{journal}{European Series in Applied and Industrial
  Mathematics (ESAIM): Mathematical Modelling and Numerical Analysis}
  \bibinfo{volume}{57}, \bibinfo{pages}{2159--2191}.
\bibitem[{Neittaanm\"{a}ki and Repin(2004)}]{neitt2004}
\bibinfo{author}{Neittaanm\"{a}ki, P.}, \bibinfo{author}{Repin, S.},
  \bibinfo{year}{2004}.
\newblock \bibinfo{title}{Reliable methods for computer simulation}.
  volume~\bibinfo{volume}{33} of \textit{\bibinfo{series}{Studies in
  Mathematics and its Applications}}.
\newblock \bibinfo{publisher}{Elsevier Science B.V., Amsterdam}.
\bibitem[{Nguyen(2024)}]{nguyen2024model}
\bibinfo{author}{Nguyen, N.C.}, \bibinfo{year}{2024}.
\newblock \bibinfo{title}{Model reduction techniques for parametrized nonlinear
  partial differential equations}.
\newblock \bibinfo{journal}{Advances in Applied Mechanics (AAMS)}
  \bibinfo{volume}{58}.
\bibitem[{Nguyen et~al.(2015)Nguyen, Anitescu, Bordas and
  Rabczuk}]{nguyen2015isogeometric}
\bibinfo{author}{Nguyen, V.P.}, \bibinfo{author}{Anitescu, C.},
  \bibinfo{author}{Bordas, S.P.}, \bibinfo{author}{Rabczuk, T.},
  \bibinfo{year}{2015}.
\newblock \bibinfo{title}{Isogeometric analysis: an overview and computer
  implementation aspects}.
\newblock \bibinfo{journal}{Mathematics and Computers in Simulation}
  \bibinfo{volume}{117}, \bibinfo{pages}{89--116}.
\bibitem[{Nguyen et~al.(2014)Nguyen, Kerfriden, Brino, Bordas and
  Bonisoli}]{nguyen2014}
\bibinfo{author}{Nguyen, V.P.}, \bibinfo{author}{Kerfriden, P.},
  \bibinfo{author}{Brino, M.}, \bibinfo{author}{Bordas, S.P.A.},
  \bibinfo{author}{Bonisoli, E.}, \bibinfo{year}{2014}.
\newblock \bibinfo{title}{Nitsche's method for two and three dimensional
  {NURBS} patch coupling}.
\newblock \bibinfo{journal}{Computational Mechanics} \bibinfo{volume}{53},
  \bibinfo{pages}{1163--1182}.
\bibitem[{Nguyen et~al.(2008)Nguyen, Rabczuk, Bordas and
  Duflot}]{nguyen2008meshless}
\bibinfo{author}{Nguyen, V.P.}, \bibinfo{author}{Rabczuk, T.},
  \bibinfo{author}{Bordas, S.P.A.}, \bibinfo{author}{Duflot, M.},
  \bibinfo{year}{2008}.
\newblock \bibinfo{title}{Meshless methods: a review and computer
  implementation aspects}.
\newblock \bibinfo{journal}{Mathematics and Computers in Simulation}
  \bibinfo{volume}{79}, \bibinfo{pages}{763--813}.
\bibitem[{Nguyen-Xuan et~al.(2008a)Nguyen-Xuan, Bordas and
  Nguyen-Dang}]{nguyen2008smooth}
\bibinfo{author}{Nguyen-Xuan, H.}, \bibinfo{author}{Bordas, S.},
  \bibinfo{author}{Nguyen-Dang, H.}, \bibinfo{year}{2008}a.
\newblock \bibinfo{title}{Smooth finite element methods: convergence, accuracy
  and properties}.
\newblock \bibinfo{journal}{International Journal for Numerical Methods in
  Engineering} \bibinfo{volume}{74}, \bibinfo{pages}{175--208}.
\bibitem[{Nguyen-Xuan et~al.(2008b)Nguyen-Xuan, Rabczuk, Bordas and
  Debongnie}]{nguyen-xuan-2008}
\bibinfo{author}{Nguyen-Xuan, H.}, \bibinfo{author}{Rabczuk, T.},
  \bibinfo{author}{Bordas, S.}, \bibinfo{author}{Debongnie, J.F.},
  \bibinfo{year}{2008}b.
\newblock \bibinfo{title}{A smoothed finite element method for plate analysis}.
\newblock \bibinfo{journal}{Computer Methods in Applied Mechanics and
  Engineering} \bibinfo{volume}{197}, \bibinfo{pages}{1184--1203}.
\bibitem[{Nochetto and Veeser(2012)}]{NochettoVeeser12}
\bibinfo{author}{Nochetto, R.H.}, \bibinfo{author}{Veeser, A.},
  \bibinfo{year}{2012}.
\newblock \bibinfo{title}{Primer of adaptive finite element methods}, in:
  \bibinfo{booktitle}{Multiscale and adaptivity: modeling, numerics and
  applications}. \bibinfo{publisher}{Springer, Heidelberg}. volume
  \bibinfo{volume}{2040} of \textit{\bibinfo{series}{Lecture Notes in Math.}},
  pp. \bibinfo{pages}{125--225}.
\bibitem[{Oden et~al.(2005)Oden, Babu{\v{s}}ka, Nobile, Feng and
  Tempone}]{oden2005}
\bibinfo{author}{Oden, J.T.}, \bibinfo{author}{Babu{\v{s}}ka, I.},
  \bibinfo{author}{Nobile, F.}, \bibinfo{author}{Feng, Y.},
  \bibinfo{author}{Tempone, R.}, \bibinfo{year}{2005}.
\newblock \bibinfo{title}{Theory and methodology for estimation and control of
  errors due to modeling, approximation, and uncertainty}.
\newblock \bibinfo{journal}{Computer Methods in Applied Mechanics and
  Engineering} \bibinfo{volume}{194}, \bibinfo{pages}{195--204}.
\bibitem[{Orkisz and Milewski(2008)}]{meshless_orkisz2008posteriori}
\bibinfo{author}{Orkisz, J.}, \bibinfo{author}{Milewski, S.},
  \bibinfo{year}{2008}.
\newblock \bibinfo{title}{A’posteriori error estimation based on higher order
  approximation in the meshless finite difference method}, in:
  \bibinfo{booktitle}{Meshfree Methods for Partial Differential Equations IV},
  \bibinfo{organization}{Springer}. pp. \bibinfo{pages}{189--213}.
\bibitem[{Panetier et~al.(2010)Panetier, Ladeveze and
  Chamoin}]{panetier2010strict}
\bibinfo{author}{Panetier, J.}, \bibinfo{author}{Ladeveze, P.},
  \bibinfo{author}{Chamoin, L.}, \bibinfo{year}{2010}.
\newblock \bibinfo{title}{Strict and effective bounds in goal-oriented error
  estimation applied to fracture mechanics problems solved with xfem}.
\newblock \bibinfo{journal}{International Journal for Numerical Methods in
  Engineering} \bibinfo{volume}{81}, \bibinfo{pages}{671--700}.
\bibitem[{Papež(2024)}]{papez2024algebraic}
\bibinfo{author}{Papež, J.}, \bibinfo{year}{2024}.
\newblock \bibinfo{title}{Algebraic error in numerical pdes and its
  estimation}.
\newblock \bibinfo{journal}{Advances in Applied Mechanics (AAMS)}
  \bibinfo{volume}{58}.
\bibitem[{Park et~al.(2003)Park, Kwon and Youn}]{meshless_park2003posteriori}
\bibinfo{author}{Park, S.H.}, \bibinfo{author}{Kwon, K.C.},
  \bibinfo{author}{Youn, S.K.}, \bibinfo{year}{2003}.
\newblock \bibinfo{title}{A posteriori error estimates and an adaptive scheme
  of least-squares meshfree method}.
\newblock \bibinfo{journal}{International Journal for Numerical Methods in
  Engineering} \bibinfo{volume}{58}, \bibinfo{pages}{1213--1250}.
\bibitem[{Pelle et~al.(1996)Pelle, Beckers and Gallimard}]{pelle1996}
\bibinfo{author}{Pelle, J.P.}, \bibinfo{author}{Beckers, P.},
  \bibinfo{author}{Gallimard, L.}, \bibinfo{year}{1996}.
\newblock \bibinfo{title}{Estimation des erreurs de discrétisation et analyses
  adaptatives, application à l'automatisation des calculs éléments finis}.
\newblock \bibinfo{note}{Institut pour la promotion des sciences de
  l'ingénieur}.
\bibitem[{Perazzo et~al.(2008)Perazzo, L{\"o}hner and
  Perez-Pozo}]{perazzo2008adaptive}
\bibinfo{author}{Perazzo, F.}, \bibinfo{author}{L{\"o}hner, R.},
  \bibinfo{author}{Perez-Pozo, L.}, \bibinfo{year}{2008}.
\newblock \bibinfo{title}{Adaptive methodology for meshless finite point
  method}.
\newblock \bibinfo{journal}{Advances in Engineering Software}
  \bibinfo{volume}{39}, \bibinfo{pages}{156--166}.
\bibitem[{Peskin(2002)}]{peskin2002}
\bibinfo{author}{Peskin, C.S.}, \bibinfo{year}{2002}.
\newblock \bibinfo{title}{The immersed boundary method}.
\newblock \bibinfo{journal}{Acta Numerica} \bibinfo{volume}{11},
  \bibinfo{pages}{479--517}.
\bibitem[{Plewa et~al.(2005)Plewa, Linde and Weirs}]{plewa2005}
\bibinfo{editor}{Plewa, T.}, \bibinfo{editor}{Linde, T.},
  \bibinfo{editor}{Weirs, V.G.} (Eds.), \bibinfo{year}{2005}.
\newblock \bibinfo{title}{Adaptive mesh refinement -- theory and applications.
  {Proceedings} of the {Chicago} workshop on adaptive mesh refinement methods,
  {Chicago}, {IL}, {USA}, {September} 3--5, 2003.}. volume~\bibinfo{volume}{41}
  of \textit{\bibinfo{series}{Lect. Notes Comput. Sci. Eng.}}
\newblock \bibinfo{publisher}{Berlin: Springer}.
\bibitem[{Prange et~al.(2012)Prange, Loehnert and Wriggers}]{prange2012error}
\bibinfo{author}{Prange, C.}, \bibinfo{author}{Loehnert, S.},
  \bibinfo{author}{Wriggers, P.}, \bibinfo{year}{2012}.
\newblock \bibinfo{title}{Error estimation for crack simulations using the
  xfem}.
\newblock \bibinfo{journal}{International Journal for Numerical Methods in
  Engineering} \bibinfo{volume}{91}, \bibinfo{pages}{1459--1474}.
\bibitem[{Prudhomme and Oden(1999)}]{prudhomme1999}
\bibinfo{author}{Prudhomme, S.}, \bibinfo{author}{Oden, J.T.},
  \bibinfo{year}{1999}.
\newblock \bibinfo{title}{On goal-oriented error estimation for elliptic
  problems: {Application} to the control of pointwise errors}.
\newblock \bibinfo{journal}{Computer Methods in Applied Mechanics and
  Engineering} \bibinfo{volume}{176}, \bibinfo{pages}{313--331}.
\bibitem[{Quarteroni and Valli(1994)}]{quarteroni1994}
\bibinfo{author}{Quarteroni, A.}, \bibinfo{author}{Valli, A.},
  \bibinfo{year}{1994}.
\newblock \bibinfo{title}{Numerical approximation of partial differential
  equations}. volume~\bibinfo{volume}{23} of \textit{\bibinfo{series}{Springer
  Series in Computational Mathematics}}.
\newblock \bibinfo{publisher}{Springer-Verlag, Berlin}.
\bibitem[{Rabczuk and Belytschko(2005)}]{meshless_rabczuk2005adaptivity}
\bibinfo{author}{Rabczuk, T.}, \bibinfo{author}{Belytschko, T.},
  \bibinfo{year}{2005}.
\newblock \bibinfo{title}{Adaptivity for structured meshfree particle methods
  in 2d and 3d}.
\newblock \bibinfo{journal}{International Journal for Numerical Methods in
  Engineering} \bibinfo{volume}{63}, \bibinfo{pages}{1559--1582}.
\bibitem[{Racz and Bui(2012)}]{meshless_integration_racz2012novel}
\bibinfo{author}{Racz, D.}, \bibinfo{author}{Bui, T.Q.}, \bibinfo{year}{2012}.
\newblock \bibinfo{title}{Novel adaptive meshfree integration techniques in
  meshless methods}.
\newblock \bibinfo{journal}{International Journal for Numerical Methods in
  Engineering} \bibinfo{volume}{90}, \bibinfo{pages}{1414--1434}.
\bibitem[{Repin(2024)}]{repin2024posteriori}
\bibinfo{author}{Repin, S.I.}, \bibinfo{year}{2024}.
\newblock \bibinfo{title}{A posteriori error identities and estimates of
  modelling errors}.
\newblock \bibinfo{journal}{Advances in Applied Mechanics (AAMS)}
  \bibinfo{volume}{58}.
\bibitem[{R{\'o}denas et~al.(2008)R{\'o}denas, Gonz{\'a}lez-Estrada,
  Taranc{\'o}n and Fuenmayor}]{rodenas2008recovery}
\bibinfo{author}{R{\'o}denas, J.J.}, \bibinfo{author}{Gonz{\'a}lez-Estrada,
  O.A.}, \bibinfo{author}{Taranc{\'o}n, J.E.}, \bibinfo{author}{Fuenmayor,
  F.J.}, \bibinfo{year}{2008}.
\newblock \bibinfo{title}{A recovery-type error estimator for the extended
  finite element method based on singular+ smooth stress field splitting}.
\newblock \bibinfo{journal}{International Journal for Numerical Methods in
  Engineering} \bibinfo{volume}{76}, \bibinfo{pages}{545--571}.
\bibitem[{Rognes and Logg(2013)}]{rognes2013}
\bibinfo{author}{Rognes, M.E.}, \bibinfo{author}{Logg, A.},
  \bibinfo{year}{2013}.
\newblock \bibinfo{title}{Automated goal-oriented error control. {I}:
  {Stationary} variational problems}.
\newblock \bibinfo{journal}{SIAM Journal on Scientific Computing}
  \bibinfo{volume}{35}, \bibinfo{pages}{c173--c193}.
\bibitem[{R{\"u}ter et~al.(2013)R{\"u}ter, Gerasimov and Stein}]{ruter2013}
\bibinfo{author}{R{\"u}ter, M.}, \bibinfo{author}{Gerasimov, T.},
  \bibinfo{author}{Stein, E.}, \bibinfo{year}{2013}.
\newblock \bibinfo{title}{Goal-oriented explicit residual-type error estimates
  in {XFEM}}.
\newblock \bibinfo{journal}{Computational Mechanics} \bibinfo{volume}{52},
  \bibinfo{pages}{361--376}.
\bibitem[{Stein(2014)}]{stein2014}
\bibinfo{author}{Stein, E.}, \bibinfo{year}{2014}.
\newblock \bibinfo{title}{History of the finite element method -- mathematics
  meets mechanics. {II}: {Mathematical} foundation of primal {FEM} for elastic
  deformations, error analysis and adaptivity}, in: \bibinfo{booktitle}{The
  history of theoretical, material and computational mechanics -- mathematics
  meets mechanics and engineering. Selected papers based on the special
  sessions at the GAMM 2010, Karlsruhe, Germany, GAMM 2011, Graz, Austria, GAMM
  2012, Darmstadt, Germany}. \bibinfo{publisher}{Berlin: Springer}, pp.
  \bibinfo{pages}{443--478}.
\bibitem[{Stevenson(2007)}]{stevenson2007}
\bibinfo{author}{Stevenson, R.}, \bibinfo{year}{2007}.
\newblock \bibinfo{title}{Optimality of a standard adaptive finite element
  method}.
\newblock \bibinfo{journal}{Foundations of Computational Mathematics}
  \bibinfo{volume}{7}, \bibinfo{pages}{245--269}.
\bibitem[{Sukumar and Tabarraei(2004)}]{sukumar2004}
\bibinfo{author}{Sukumar, N.}, \bibinfo{author}{Tabarraei, A.},
  \bibinfo{year}{2004}.
\newblock \bibinfo{title}{Conforming polygonal finite elements}.
\newblock \bibinfo{journal}{International Journal for Numerical Methods in
  Engineering} \bibinfo{volume}{61}, \bibinfo{pages}{2045--2066}.
\bibitem[{Szab{\'o} and Babu{\v{s}}ka(2021)}]{szabo2021}
\bibinfo{author}{Szab{\'o}, B.}, \bibinfo{author}{Babu{\v{s}}ka, I.},
  \bibinfo{year}{2021}.
\newblock \bibinfo{title}{Finite element analysis. {Formulation}, verification
  and validation}.
\newblock Wiley Series Comput. Mech.. \bibinfo{edition}{2nd updated and revised
  edition} ed., \bibinfo{publisher}{Hoboken, NJ: John Wiley \& Sons}.
\newblock \DOIprefix\doi{10.1002/9781119426479}.
\bibitem[{Veeser and Verf\"{u}rth(2009)}]{veeser2009}
\bibinfo{author}{Veeser, A.}, \bibinfo{author}{Verf\"{u}rth, R.},
  \bibinfo{year}{2009}.
\newblock \bibinfo{title}{Explicit upper bounds for dual norms of residuals}.
\newblock \bibinfo{journal}{SIAM Journal on Numerical Analysis}
  \bibinfo{volume}{47}, \bibinfo{pages}{2387--2405}.
\bibitem[{Veeser and Verf\"{u}rth(2012)}]{veeser2012}
\bibinfo{author}{Veeser, A.}, \bibinfo{author}{Verf\"{u}rth, R.},
  \bibinfo{year}{2012}.
\newblock \bibinfo{title}{Poincar\'{e} constants for finite element stars}.
\newblock \bibinfo{journal}{IMA Journal of Numerical Analysis}
  \bibinfo{volume}{32}, \bibinfo{pages}{30--47}.
\bibitem[{Verf\"{u}rth(1999)}]{verfurth1999}
\bibinfo{author}{Verf\"{u}rth, R.}, \bibinfo{year}{1999}.
\newblock \bibinfo{title}{Error estimates for some quasi-interpolation
  operators}.
\newblock \bibinfo{journal}{M2AN. Mathematical Modelling and Numerical
  Analysis} \bibinfo{volume}{33}, \bibinfo{pages}{695--713}.
\bibitem[{Verf\"{u}rth(2009)}]{verfurth2009}
\bibinfo{author}{Verf\"{u}rth, R.}, \bibinfo{year}{2009}.
\newblock \bibinfo{title}{A note on constant-free a posteriori error
  estimates}.
\newblock \bibinfo{journal}{SIAM Journal on Numerical Analysis}
  \bibinfo{volume}{47}, \bibinfo{pages}{3180--3194}.
\bibitem[{Verf{\"u}rth(2013)}]{verfurth2013}
\bibinfo{author}{Verf{\"u}rth, R.}, \bibinfo{year}{2013}.
\newblock \bibinfo{title}{A posteriori error estimation techniques for finite
  element methods}.
\newblock Numer. Math. Sci. Comput., \bibinfo{publisher}{Oxford: Oxford
  University Press}.
\bibitem[{Videla et~al.(2019)Videla, Anitescu, Khajah, Bordas and
  Atroshchenko}]{videla2019h}
\bibinfo{author}{Videla, J.}, \bibinfo{author}{Anitescu, C.},
  \bibinfo{author}{Khajah, T.}, \bibinfo{author}{Bordas, S.P.},
  \bibinfo{author}{Atroshchenko, E.}, \bibinfo{year}{2019}.
\newblock \bibinfo{title}{h-and p-adaptivity driven by recovery and
  residual-based error estimators for pht-splines applied to time-harmonic
  acoustics}.
\newblock \bibinfo{journal}{Computers \& Mathematics with Applications}
  \bibinfo{volume}{77}, \bibinfo{pages}{2369--2395}.
\bibitem[{Videla et~al.(2024)Videla, Shaaban and
  Atroshchenko}]{IGA_GIFT_videla2024shape}
\bibinfo{author}{Videla, J.}, \bibinfo{author}{Shaaban, A.M.},
  \bibinfo{author}{Atroshchenko, E.}, \bibinfo{year}{2024}.
\newblock \bibinfo{title}{Shape optimization with adaptive geometry independent
  field approximation (gift) in 3d time-harmonic acoustics}.
\newblock \bibinfo{journal}{Journal of Sound and Vibration} ,
  \bibinfo{pages}{118299}.
\bibitem[{Vohral{\'{\i}}k(2007)}]{Vohralik2007SIAMJNA}
\bibinfo{author}{Vohral{\'{\i}}k, M.}, \bibinfo{year}{2007}.
\newblock \bibinfo{title}{A posteriori error estimates for lowest-order mixed
  finite element discretizations of convection-diffusion-reaction equations}.
\newblock \bibinfo{journal}{SIAM Journal on Numerical Analysis}
  \bibinfo{volume}{45}, \bibinfo{pages}{1570--1599}.
\bibitem[{Yu et~al.(2018)Yu, Anitescu, Tomar, Bordas and
  Kerfriden}]{yu2018adaptive}
\bibinfo{author}{Yu, P.}, \bibinfo{author}{Anitescu, C.},
  \bibinfo{author}{Tomar, S.}, \bibinfo{author}{Bordas, S.P.A.},
  \bibinfo{author}{Kerfriden, P.}, \bibinfo{year}{2018}.
\newblock \bibinfo{title}{Adaptive isogeometric analysis for plate vibrations:
  an efficient approach of local refinement based on hierarchical a posteriori
  error estimation}.
\newblock \bibinfo{journal}{Computer Methods in Applied Mechanics and
  Engineering} \bibinfo{volume}{342}, \bibinfo{pages}{251--286}.
\bibitem[{Zienkiewicz and Zhu(1987)}]{zienkiewicz1987}
\bibinfo{author}{Zienkiewicz, O.C.}, \bibinfo{author}{Zhu, J.Z.},
  \bibinfo{year}{1987}.
\newblock \bibinfo{title}{A simple error estimator and adaptive procedure for
  practical engineering analysis}.
\newblock \bibinfo{journal}{International Journal for Numerical Methods in
  Engineering} \bibinfo{volume}{24}, \bibinfo{pages}{337--357}.

\end{thebibliography}

\end{document}